\documentclass[10pt,reqno]{amsart}
\usepackage{amssymb}
\usepackage{amsmath}
\usepackage{bm}
\usepackage{amsthm}
\usepackage{amsbsy}
\usepackage{psfrag}
\usepackage{graphics}
\usepackage{graphicx}
\usepackage{pgfplots}
\usepackage{relsize}
\usepackage{color}
\usepackage{xcolor}
\usepackage{pdfpages}
\usepackage{url}
\usepackage{hyperref}
\usepackage[square,sort&compress,comma,numbers]{natbib}
\usepackage{graphicx}
\usepackage{caption}
\usepackage{subcaption}
\hypersetup{
	colorlinks=true,
	linkcolor=blue,
	filecolor=magenta,      
	urlcolor=cyan,
}

\theoremstyle{plain}
\theoremstyle{theorem}
\newtheorem{theorem}{Theorem}[section]
\newtheorem{proposition}[theorem]{Proposition}

\newtheorem{lemma}[theorem]{Lemma}

\newtheorem{example}[theorem]{Example}
\newcommand{\norm}[1]{\left\Vert #1 \right\Vert}
\newtheorem{remark}[theorem]{Remark}

\begin{document}

%
%
%
%
%

\title[FEM for Dirichlet Boundary Optimal Control  Problem]
{Finite Element Analysis of the Dirichlet Boundary Control Problem Governed by Linear Parabolic Equation}

\author[T. Gudi]{ Thirupathi Gudi}
\address{Department of Mathematics, Indian Institute of Science, Bangalore - 560012, India}
\email{gudi@iisc.ac.in}

\author[G. Mallik]{ Gouranga Mallik}
\address{Department of Mathematics, Indian Institute of Science, Bangalore - 560012, India}
\email{gourangam@iisc.ac.in}

\author[R. C. Sau]{Ramesh Ch. Sau}
\address{Department of Mathematics, Indian Institute of Science, Bangalore - 560012, India}
\email{rameshsau@iisc.ac.in}

\date{}

\begin{abstract}
  A finite element analysis of a Dirichlet boundary control problem governed by the linear parabolic equation is presented in this article. The Dirichlet control is considered in a closed and convex subset of the energy space $H^1(\Omega \times(0,T)).$ We prove well-posedness and discuss some regularity results for the control problem. We derive the optimality system for the optimal control problem. The first order necessary optimality condition results in a simplified Signorini type problem for control variable. The space discretization of the state variable is done using conforming finite elements, whereas the time discretization is based on discontinuous Galerkin methods. To discretize the control we use the conforming prismatic Lagrange finite elements. We derive an optimal order of convergence of error in control, state and adjoint state. The theoretical results are corroborated by some numerical tests.	
\end{abstract}
\keywords{ PDE-constrained optimization; Control-constraints;  Finite element
	method; Error bounds; Evolution equation}

\subjclass{65N30; 65N15; 65N12; 65K10}

\maketitle
\def \R{{{\Bbb R}}}
\def \P{{{\rm {\!}\cal P}}}
\def \Z{{{\rm {\!}\cal Z}}}
\def \Q{{{\rm {\!}\cal Q}}}
\allowdisplaybreaks
\def\d{\displaystyle}
\def\R{\mathbb{R}}
\def\I{\mathbb{I}}
\def\cA{\mathcal{A}}
\def\dx{{\rm~dx}}
\def\dy{{\rm~dy}}
\def\ds{{\rm~ds}}
\def\y{\bar{q}}
\def\u{\bar{u}}
\def\V{\mathbf{V}}
\def\w{\bar{w}}
\def\J{\mathbf{J}}
\def\f{\mathbf{f}}
\def\g{\mathbf{g}}
\def\div{{\rm div~}}
\def\v{\mathbf{v}}
\def\z{\mathbf{z}}
\def\b{\mathbf{b}}
\def\cB{\mathcal{B}}
\def\p{\partial}
\def\O{\Omega}
\def\s{\sigma}
\def\bbP{\mathbb{P}}
\def\vb{{\vec b}}
\def\bn{{\bf n}}
\def\cV{\mathcal{V}}
\def\V{\mathbf{V}}
\def\cR{\mathcal{R}}
\def\cI{\mathcal{I}}
\def\cM{\mathcal{M}}
\def\cT{\mathcal{T}}
\def\cE{\mathcal{E}}
\def\bbE{\mathbb{E}}
\def\ssT{{\scriptscriptstyle T}}
\def\HT{{H^2(\O,\cT_h)}}
\def\mean#1{\left\{\hskip -5pt\left\{#1\right\}\hskip -5pt\right\}}
\def\jump#1{\left[\hskip -3.5pt\left[#1\right]\hskip -3.5pt\right]}
\def\smean#1{\{\hskip -3pt\{#1\}\hskip -3pt\}}
\def\sjump#1{[\hskip -1.5pt[#1]\hskip -1.5pt]}
\def\jumptwo{\jump{\frac{\p^2 u_h}{\p n^2}}}

\section{Introduction}\label{sec:Intro:Parabolic}
The study of optimal control problem govern by partial differential equations (PDEs) is a significant area of research in applied mathematics. The optimal control problem consists of finding a control variable that minimizes a cost functional subject to a PDE. 
Due to the importance in applications, several numerical methods have been proposed to approximate the solutions. The finite element approximation of the optimal control problem started with the work of Falk \cite{Falk:1973:Control} and Geveci \cite{Gevec:1979:Control}. 
A control can act in the interior of a domain, in this case, we call distributed, or on the boundary of a domain, we call boundary (Neumann or Dirichlet) control problem. We refer to \cite{Hinze:2011:Control,meyerrosch:2004:Optiaml,Sudipto:2015:IP,Ramesh:2018} for distributed control related problem, to \cite{CasasM:2008:Neum,Casas_Dhamo_2012,Sudipto:2015:IP,Ramesh:2018} for the Neumann boundary control problem, and to \cite{Hinze:2005:Control} for a variational discretization approach. The Dirichlet boundary control problem has been studied in \cite{CMR:2009:Dirich,CasasRaymond:2006:Dirich}.

The Dirichlet boundary control problems are essential in the application areas, and various approaches are proposed in the literature for the same.
One such is to seek control from $L^2(\Gamma)$-space (see \cite{CasasRaymond:2006:Dirich}). In this case the state equation  has to be understood in a ultra weak sense, since the Dirichlet boundary data is only in $L^2(\Gamma)$. This ultra-weak formulation is easy to implement and typically yields optimal controls with low regularity. Especially, when the problem is posed on a  polygonal domain, the control exhibits layer behaviour at the corner points. This is because it is determined by the normal derivative of the adjoint state.
Another approach is to choose the control from the the energy space $H^{1/2}(\Gamma)$ (see \cite{Steinbach:2014:Dirichlet}). With the help of a harmonic extension of the given boundary data, the Steklov-Poincar\'e operator was employed in \cite{Steinbach:2014:Dirichlet} to determine the cost functional. By employing harmonic extension of the Dirichlet data, the Steklov-Poincar\'e operator turns Dirichlet data into Neumann data; nevertheless, numerical implementation of this sort of abstract operator might be challenging. In paper \cite{CMR:2009:Dirich}, Dirichlet control problem is transformed into a Robin boundary control problem through penalization. In \cite{Gudi:2017:DiriControl}, the authors consider unconstrained Dirichlet boundary control where  the control in $H^{1/2}(\Gamma)$ is realized by a harmonic extension in $H^1(\Omega)$ which enables to consider cost functional in energy form. In this approach the authors choose the control from the energy space $H^1(\Omega)$ so that  they do not need Steklov-Poincar\'e operator and hence this method is computationally very efficient.  We refer \cite{Gudi_Ramesh:ESIAM_COCV} for an improved analysis of constrained Dirichlet boundary control.


In \cite{Kunisch_Vexler:2007}, a semi-smooth Newton method has been used to solve Dirichlet boundary control problem for parabolic PDE.  The article \cite{Arada-Raymond:2002,Belgacem-bernadi:2011} investigates the Robin-type boundary conditions for parabolic Dirichlet boundary control problems using Robin penalization method. In this paper, we consider the following parabolic Dirichlet boundary control problem of tracking type, which may be regarded as prototype problem(based on energy approach) to study Dirichlet boundary control for time-dependent PDEs.	
\begin{align}\label{continuous_cost_functional}
	\text{min} ~J(u,q)= \frac{1}{2}\norm{u-u_d}_{L^2(I;L^2(\Omega))}^2 + \frac{\lambda}{2}|q|^2_{1,\Omega\times (0,T)},
\end{align}
subject to PDE,
\begin{subequations}\label{Eq:state_conti_intro}
	\begin{align}
		\partial_tu-\Delta u&=f \quad \text{in}\;\Omega \times (0,T),\label{parabolic_pde} \\
		u&=q \quad \text{on}\; \;\partial\Omega \times (0,T),\\
		u(x,0)&=u_0(x) \quad \text{in}\; \Omega,
	\end{align}
\end{subequations} 
with the control constraints  $$ q_a\leq q(x,t)\leq q_b \quad \text{on}\;\partial\Omega \times (0,T).$$

\noindent
The detailed description of the above problem will be discussed in the Section \ref{sec:ModelProblem:para}. To the authors' knowledge, this is the first work to address the energy approach for solving the Dirichlet boundary control problem governed by linear parabolic equation. We prove existence and uniqueness of the solution of the problem \eqref{continuous_cost_functional}, and discuss about the regularity of the optimal control, which is of particular interest. Also, one of the goal of this article is to consider the discretization of the optimality system based on a finite element approximation of the state, adjoint state and the control variable. To discretize the state and adjoint state equation, we use discontinuous finite elements for time discretization, and $H^1$-conforming finite elements for spatial discretization.  In \cite{Vexler_et.al:2007} this type of discretization is shown to allow for a natural translation of the optimality conditions from the continuous to the discrete level. This gives rise to exact computation of the derivatives required in the optimization algorithms on the discrete level.  
Since, we sought control from a closed convex subset of  $H^1(\Omega\times (0,T))$, the optimal control satisfies a simplified Signorini problem in three dimensional domain $\Omega\times (0,T)$. To discretize the control we use the conforming prismatic Lagrange finite elements on three dimensional domain. For a period of fifteen years, the lack of order of convergence for the Signorini solution had been a common difficulty. Hild \cite{Hild_15_signorini},  derive the optimal order of convergence under minimal assumptions. Using the ideas from \cite{Hild_15_signorini}, we derive the optimal order convergence for the error in control.

The following is a breakdown of the rest of the article. The investigated Dirichlet boundary control problem, the primal boundary value problem, the reduced cost functional, and the related adjoint boundary value problem are all described in Section \ref{sec:ModelProblem:para}. The minimizer of the reduced cost functional is characterized as the unique solution of a variational inequality of the first kind. Section \ref{discretization} discusses the finite element discretization of the variational inequality, as well as finite element approximations of both the primal and adjoint boundary value problems, and related error estimates.  Some numerical results are finally given in Section \ref{sec5:para}.

\section{Continuous Control Problem}\label{sec:ModelProblem:para}
In this section, we first introduce some notations and then discuss the  mathematical formulation of the optimal control problem. Furthermore, we prove theoretical results on existence and uniqueness.
\subsection{Notations and Preliminary}\label{notation} Let $\Omega$ be a bounded convex polygonal domain in $\mathbb{R}^2.$ The inner-product in $L^2(\Omega)$ space is denoted by $(u,v):=\int_{\Omega} uv$ with the norm $\norm{v}_{0,\Omega}:=(v,v)^{1/2}$. Also, we denote the norm in the $H^k(\Omega)$ space by $\norm{v}_{k,\Omega}$ for $k\geq 1$. Let 
$I:=(0,T)$ be a time interval. We consider the Sobolev space $L^2(I;L^2(\Omega))$ equipped with the inner-product $$(u,v)_I=\int_{0}^{T}\int_{\Omega}uv$$ and the norm  $\norm{u}_{I}:=(u,u)_{I}^{1/2}$.  Let $0<s<1$.  Define the fractional Sobolev space
$$H^{s}(\Omega):=\{ u \in L^2(\Omega)\,|\; \int_{\Omega}\int_{\Omega} \frac{(u(x)-u(y))^2}{(x-y)^{2s+2}}\dx\dy< \infty \},$$
equipped with the norm
$$\norm{u}_{s, \Omega}:=\bigg{(}\norm{u}_{0,\Omega}^2+ \int_{\Omega}\int_{\Omega} \frac{(u(x)-u(y))^2}{(x-y)^{2s+2}}\dx\dy\bigg{)}^{\frac{1}{2}}.$$
The fractional Sobolev space involving time is defined by
$$L^2(I;H^s(\Omega)):=\{u:I\rightarrow H^s(\Omega)\; \text{measurable}\;|\; \int_{I}\norm{u(t)}^2_{s,\Omega}dt\;<\infty\}.$$

Prior to discussing the optimal control problem, we describe the problem setup for the linear parabolic problem defined in \eqref{Eq:state_conti_intro} as follows:
\begin{subequations}\label{Eq:state_conti}
	\begin{align}
		\partial_tu-\Delta u&=f \quad \text{in}\;\Omega \times (0,T), \\
		u&=q \quad \text{on}\; \;\partial\Omega \times (0,T),\\
		u(x,0)&=u_0(x) \quad \text{in}\; \Omega,
	\end{align}
\end{subequations}
where $\Omega\subset \mathbb{R}^2,$ be a bounded polygonal domain and $(0,T)$ be a time interval. Denote $\Gamma_{D}:=\Omega\times\{0,T\}$ and $\Gamma_C:=\partial\Omega \times (0,T)$.
We take the interior force $f\in L^2(I;L^2(\Omega))$ and the initial data $u_0\in H^1(\Omega)$. The Dirichlet data function $q$ is the control variable and it is chosen from the following admissible space:
\begin{align}\label{control_space}
	Q=\{q\in H^1(\Omega\times I )\; |\;q(x,t)=0\; \text{on}\;\Gamma_{D}\}.
\end{align}
\noindent
For the time interval $I$,  define the test and trial space
\begin{align*}
	X:=\{w\;|\;w\in L^2(I;H_0^1(\Omega))\; \text{and}\; \partial_tw \in L^2(I;H^{-1}(\Omega))\}.
\end{align*}
For a given control $q\in Q,$ the weak formulation of \eqref{Eq:state_conti} is to find $w(q)\in X,$ with $u(x,0)=u_0(x),$ such that
\begin{align}\label{Eq:state_weak_form}
	(\partial_t w(q),v)_I+(\nabla w(q), \nabla v)_I=&(f,v)_I-(\nabla q, \nabla v)_I-(\partial_t q,v)_I \quad \forall v \in X,
\end{align}
and set $u(q):=w(q)+q\in X+Q$ to be the weak solution of \eqref{Eq:state_conti}.

We recall the following result on existence, uniqueness (see \cite{Lions:MagenesII,Book:Wolka})
and regularity (see \cite{Book:Evans,Lions:MagenesII}) for the state equation.
\begin{proposition}
	For a given control $q\in Q$, $f\in L^2(I;L^2(\Omega))$ and $u_0\in H^1(\Omega),$ there exists a unique weak solution $u(q):=w(q)+q\in X+Q$ of the problem \eqref{Eq:state_conti}. Moreover, the solution exhibits the
	improved regularity $u(q)\in Z+Q$ and satisfies the stability estimate
	\begin{align}
		\norm{\partial_t u(q)}_I+\norm{\nabla u(q)}_I\leq \norm{f}_I+\norm{\partial_t q}_I+\norm{\nabla q}_I+\norm{u_0}_{1,\Omega},
	\end{align}	
	where the space $Z$ is defined by 
	\begin{align*}
		Z:=\{w \; |\; w\in L^2(I;H_0^1(\Omega))\; \text{and}\; \partial_t w \in L^2(I;L^2(\Omega))\}.
	\end{align*}
\end{proposition}

\subsection{The Dirichlet Control Problem} Consider the cost functional $J$ defined in \eqref{continuous_cost_functional}. The model control
problem consists of finding $(\bar{u}(\bar{q}),\bar{q})\in (X+Q) \times Q_{ad},$ such that
\begin{align}\label{eq:min-J1:para}
	J(\bar{u}(\bar{q}),\bar{q})=\min_{(u(q),q)\in (X+Q)\times Q_{ad}} J(u(q),q),
\end{align}
subject to $u(q)=w(q)+q$  satisfying \eqref{Eq:state_weak_form}. The admissible constrained set of control reads
\begin{align}\label{constrained_set}
	Q_{ad}=\{q\in Q\; |\;q_a\leq q(x,t)\leq q_b \quad \text{on}\;\Gamma_{C} \},
\end{align}
where $q_a, q_b \in \mathbb{R}$ and for consistency assume $q_a\leq 0$ and $q_b\geq0,$ so that the admissible set $Q_{ad}$ is nonempty.

\begin{theorem}[Existence and uniqueness of control]
	There exists a unique solution of the control problem \eqref{eq:min-J1:para}.
\end{theorem}
\begin{proof}
	The cost functional $J$ is non negative. Set  $$\alpha=\inf_{(u(q),q)\in (X+Q)\times Q_{ad}} J(u(q),q).$$ Then there exists a minimizing sequence $(u_n(q_n),q_n)$ such that, $J(u_n(q_n),q_n)$ converges to $\alpha$. For the notational simplicity we denote $u_n(q_n)$ by $u_n$. Since the sequence $J(u_n,q_n)$ is convergent, the components $\|u_n-u_d\|_I$ and  $|
	q_n|_{1,\Omega\times I}$ are  bounded. As $q_n\in Q_{ad}$, using the Poincar\'e inequality we conclude that the sequence $q_n$ is bounded in $Q$. Then there
	exists a subsequence of $q_n$, still indexed by $n$ to simplify the notation, and a function $\bar{q}$,
	such that $q_n$ converges to $\bar{q}$ weakly in $Q$. It is clear that the set $Q_{ad}$ is closed and convex, the function $\bar{q}\in Q_{ad}$. An \textit{a priori} estimate of the problem \eqref{Eq:state_weak_form} yields
	\begin{align}
		\norm{\partial_t w_n}_I+\norm{\nabla w_n}_I\leq C (\norm{f}_I+|q_n|_{1,\Omega\times I}+\norm{u_0}_{1,\Omega}),
	\end{align}
	where $u_n=w_n+q_n$. Using this \textit{a priori} estimate and the boundedness of the sequence $q_n$ we conclude that the sequence $w_n$ is bounded in $H_0^1(\Omega \times I)$. We extract a subsequence, name it $w_n$ and it converges weakly to $\bar{w}$ in $H_0^1(\Omega \times I)$. Now we  show that $\bar{w}$ is the corresponding candidate for the control $\bar{q}$. 
	From \eqref{Eq:state_weak_form}, we have 
	\begin{align}\label{eq:sequential}
		(\partial_t w_n,v)_I+(\nabla w_n, \nabla v)_I=&(f,v)_I-(\nabla q_n, \nabla v)_I-(\partial_t q_n,v)_I \quad \forall v \in L^2(I;H_0^1(\Omega)). 
	\end{align}
	\noindent
	Using the above weak convergences in \eqref{eq:sequential}, we obtain 
	\begin{align*}
		(\partial_t \bar{w},v)_I+(\nabla \bar{w}, \nabla v)_I=&(f,v)_I-(\nabla \bar{q}, \nabla v)_I-(\partial_t \bar{q},v)_I \quad \forall v \in L^2(I;H_0^1(\Omega)). 
	\end{align*}
	Hence, $\bar{u}=\bar{w}+\bar{q}$ is the corresponding state for the control $\bar{q}$. The sequence $q_n$ converges to $\bar{q}$ weakly in $Q$. Therefore, it converges strongly in $L^2(\Omega \times I)$. Now, $w_n$ converges to $\bar{w}$ weakly in $H_0^1(\Omega \times I)$. Thus, it converges strongly in $L^2(\Omega \times I)$. So, $w_n+q_n=u_n$ converges strongly to $\bar{w}+\bar{q}=\bar{u}$ in $L^2(\Omega \times I)$. Using the weak lower semi continuity of the norm, we obtain $ |\bar{q}|_{1,\Omega\times I}\leq \liminf_{n\rightarrow \infty} |q_n|_{1,\Omega\times I}$. Hence, we have 
	$$J(\bar{u},\bar{q})\leq  \lim_{n\rightarrow \infty}\frac{1}{2} \|u_n-u_d\|_{I}^2+\frac{\lambda}{2}\liminf_{n\rightarrow \infty} |q_n|_{1,\Omega\times I}^2=\alpha.$$
	This proves the existence of a control $\bar{q}$ such that $J(\bar{u},\bar{q})=\alpha$.  The uniqueness of the solution follows from the strict convexity of the cost functional. 
\end{proof}

\begin{proposition}[Continuous Optimality System]\label{prop:Wellposed-C:prop} The state, adjoint state, and control $(\bar{u}(\bar{q}),\bar{\phi}(\bar{q}),\bar{q})\in (X+Q) \times X \times Q_{ad}$ satisfy the optimality system 
	\begin{subequations}
		\begin{align}
			\bar{u}(\bar{q})=&\bar{w}(\bar{q})+\bar{q},\quad \bar{w}(\bar{q})\in X,\label{eq:state_0a:para}\\
			(\partial_t\bar{w}(\bar{q}),v)_I+( \nabla \bar{w}(\bar{q}) , \nabla v)_I=&(f,v)_I -(\nabla \bar{q} , \nabla v)_I-(\partial_t\bar{q},v)_I \quad \forall v \in X,\label{eq:state_a:para}\\
			-(\partial_t \bar{\phi}(\bar{q}),v)_I+(\nabla \bar{\phi}(\bar{q}), \nabla v)_I=&(\bar{u}(\bar{q})-u_d,v)_I \quad \forall v \in L^2(I;H_0^1(\Omega)),\label{eq:adjstate_a:para}\\
			\lambda(\partial_t \bar{q},\partial_t (p-\bar{q}))_{I}+\lambda (\nabla \bar{q}, \nabla (p-\bar{q}))_I\geq&  (\bar{\phi}(\bar{q}),\partial_t(p-\bar{q}))_I+(\nabla \bar{\phi}(\bar{q}), \nabla (p-\bar{q}))_I\nonumber\\&-(\bar{u}(\bar{q})-u_d,p-\bar{q})_I \quad \forall p\in Q_{ad},\label{conti_VI:para}  
		\end{align}
	\end{subequations}
with $\bar{u}(x,0)=u_0(x)$ and $\bar{\phi}(\bar{q})(x,T)=0.$
\end{proposition}
\begin{proof}
	Consider the Lagrangian function
	\begin{align}
		\mathcal{L}(w,\phi,q)=\frac{1}{2}\norm{w+q-u_d}_I^2 + \frac{\lambda}{2}\norm{\nabla q}_I^2+\frac{\lambda}{2}\norm{\partial_t q}_I^2-\int_{I}\partial_t w \; \phi \nonumber\\-\int_{I} \nabla w {\cdot} \nabla \phi +\int_{I}f\phi -\int_{I}\partial_t q\; \phi -\int_{I} \nabla q {\cdot} \nabla \phi.
	\end{align}
	Differentiating $\mathcal{L}$ with respect to $\phi$ and $w$ at $(\bar{w}(\bar{q}),\bar{\phi}(\bar{q}),\bar{q})$ and equating to $0$, we obtain state \eqref{eq:state_a:para} and adjoint state \eqref{eq:adjstate_a:para}  respectively. 
	Now, differentiating $\mathcal{L}$ with respect to $q$ at $(\bar{w}(\bar{q}),\bar{\phi}(\bar{q}),\bar{q})$ in the direction $(p-\bar{q}),$ we get the following:
	\begin{align}
		\mathcal{D}_q\mathcal{L}(\bar{w}(\bar{q}),\bar{\phi}(\bar{q}),\bar{q})(p-\bar{q})=\int_{I}(\bar{w}(\bar{q})+\bar{q}-u_d)(p-\bar{q})+\lambda\int_{I}\nabla \bar{q}{\cdot}\nabla (p-\bar{q})\nonumber\\+\lambda\int_{I}\partial_t \bar{q}\;\partial_t (p-\bar{q})-\int_{I}\partial_t (p-\bar{q})\; \bar{\phi}(\bar{q}) -\int_{I} \nabla (p-\bar{q}) {\cdot} \nabla \bar{\phi}(\bar{q}).
	\end{align}
	The first order necessary optimality condition $\mathcal{D}_q\mathcal{L}(\bar{w}(\bar{q}),\bar{\phi}(\bar{q}),\bar{q})(p-\bar{q})\geq 0 \;\text{for all} \; p \in Q_{ad}$ yields the inequality \eqref{conti_VI:para}.
\end{proof}

It is easy to prove from \eqref{conti_VI:para} that the optimal control $\bar{q}$ solves the following Signorini problem:
\begin{subequations}\label{signorini_L2_f:para}
	\begin{align}
		-\lambda(\partial_{tt}+\Delta) \bar{q}&=0\quad \text{in}\quad \Omega \times (0,T),\\
		q_a\leq \bar{q} &\leq q_b \quad \text{on}\quad \Gamma_{C},\\
		\bar{q}&=0\quad \text{in}\;\Gamma_{D},
	\end{align} 
	and further the following holds for almost every $(x,t)\in \Gamma_{C}$:
	\begin{align}
		\text{if}\;  q_a<\bar{q}(x,t) <q_b \quad \text{then} \quad \big(\lambda\frac{\partial \bar{q}}{\partial n}-\frac{\partial \bar{\phi}(\bar{q})}{\partial
			n}\big)(x,t)&=0,\\
		\text{if}\; q_a\leq \bar{q}(x,t)<q_b \quad \text{then} \quad \big(\lambda\frac{\partial \bar{q}}{\partial n}-\frac{\partial \bar{\phi}(\bar{q})}{\partial
			n}\big)(x,t)&\geq 0,\\
		\text{if}\;  q_a<\bar{q}(x,t)\leq q_b \quad \text{then} \quad \big(\lambda\frac{\partial \bar{q}}{\partial n}-\frac{\partial \bar{\phi}(\bar{q})}{\partial
			n}\big)(x,t)&\leq 0.
	\end{align}
\end{subequations}

\begin{remark}[Regularity of Signorini problem]\label{Rem_reg_signorini:para}
	The numerical analysis of any finite element method applied to the Signorini problem \eqref{signorini_L2_f:para} requires
	the knowledge of the regularity of the solution $\bar{q}$. The Signorini condition may generate some singular behavior at the neighborhood of $\Gamma_C$, see \cite{Moussaoui_et_al}. There are many factors that affect the regularity of the solution to the Signorini problem. Some of those factors are the regularity of the data, the mixed boundary conditions
	(e.g., Neumann-Dirichlet transitions), the corners in polygonal domains and the Signorini condition which generates singularities at contact-noncontact transition points. Let $\mathbf{p}$ be a contact-noncontact transition point in the interior of  $\Gamma_{C},$ then the solution of Signorini problem \eqref{signorini_L2_f:para} $\bar{q} \in H^{\tau}(V_\mathbf{p})$ with $\tau < \frac{5}{2}$ and $V_\mathbf{p}$ be an open neighborhood of $\mathbf{p}$ (see \cite[subsection 2.3]{Belgacem:2013}, \cite[section 2]{Belhachmi:2003} and \cite{Moussaoui_et_al}). Let $\mathbf{p}\in \bar{\Gamma}_{C}\cap \bar{\Gamma}_{D}$ and $V_\mathbf{p}$ be a neighborhood of $\mathbf{p}$ in $\Omega$ such that $\bar{q}$ vanishes on $V_\mathbf{p}\cap\Gamma_{C}$ then the elliptic regularity theory on convex domain yields $\bar{q}\in H^{2}(V_\mathbf{p})$(see \cite[subsection 2.3]{Belgacem:2013} and \cite{Orlt:1995}). Now if $\bar{q}$ does not vanish on $V_\mathbf{p}\cap\Gamma_{C},$ then $\mathbf{p}$ be a contact-noncontact type transition point and hence $\bar{q} \in H^{\tau}(V_\mathbf{p})$ with $\tau < 5/2$ (see \cite[subsection 2.3]{Belgacem:2013} and \cite{Orlt:1995}). The best we can expect is to obtain $\bar{q} \in H^{\tau}(V_{\Gamma_{C}})$ with $\tau \leq 2$ and $V_{\Gamma_{C}}$ is an open neighbourhood of $\Gamma_{C}$ (see \cite{Moussaoui_et_al,Belgacem:2013}).
\end{remark}

\section{Discretization and error analysis}\label{discretization}
In this section, we consider finite element discretization of the optimal control problem \eqref{eq:min-J1:para}. Also for the error analysis we assume that the solutions  $\bar{w}(\bar{q}),\bar{\phi}(\bar{q})\in H^1(I,H^1_0(\Omega))\cap L^2(I,H^{\tau}(\Omega))$ and $\bar{q}\in H_D^{\tau}(\Omega \times I):=\{p \in H^{\tau}(\Omega \times I): p=0 \; \text{on} \;\Gamma_D\}  $ with $3/2<\tau\leq 2$.
\subsection{Discretization in time and space}\label{subsec:Discretization} 
In this subsection, we first discretize the time and then discretize the space.

\noindent{\bf Semi-discretization in time}.
Let $\bar{I}= \{0\} \cup I_1 \cup I_2 \cup ... \cup I_M$ be a partition of  $\bar{I}= [0,T]$ with subintervals $I_m = (t_{m-1},t_m]$ of length $k_m:=t_m-t_{m-1}$ and time points
\begin{align}\label{Time_partition}
	0 = t_0 < t_1 <... < t_{M-1} < t_M = T.
\end{align}
The time discretization parameter is defined by $k=\max_{1\leq m\leq M} k_m$.
The semidiscrete test and trial space is defined by
\begin{align*}
	X_k^0:=\{v_k\in L^2(I,H_0^1(\Omega))\;|\;v_k|_{I_m}\in \mathcal{P}_0(I_m,H_0^1(\Omega))\; \text{for}\;m=1,2,3,...,M\},
\end{align*}
where, $\mathcal{P}_0(I_m,H_0^1(\Omega))$ denotes the space of constant polynomials defined on $I_m$ with
values in $H_0^1(\Omega)$. For $u,v\in X_k^0$, we use the notation
$$(u,v)_{I_m}=(u,v)_{L^2(I_m,L^2(\Omega))}\quad \text{and} \quad \norm{u}_{I_m}=\norm{u}_{L^2(I_m,L^2(\Omega))}.$$ 
For $v_k \in X_k^0,$ we define the following notations:
$$v^{+}_{k,m}:=\lim_{t\rightarrow 0^{+}}v_k(t_m+t),\quad v^{-}_{k,m}:=\lim_{t\rightarrow 0^{+}}v_k(t_m-t)=v_k(t_m),\quad \jump{v_k}_m:=v^{+}_{k,m}-v^{-}_{k,m}.$$
\noindent
Define the bilinear form $B:X_k^0\times X_k^0\to\mathbb{R}$, 
\begin{align}\label{defn:B}
	B(w_k,v_k):=\sum_{m=1}^{M}\big(\partial_t w_k,v_k\big)_{I_m}+\big(\nabla w_k , \nabla v_k\big)_{I}+\sum_{m=1}^{M-1} \big(\jump{w_k}_m,v^{+}_{k,m}\big)+\big(w^{+}_{k,0},v^{+}_{k,0}\big),
\end{align}
for $w_k,v_k\in X_k^0$.
The semi-discrete weak formulation of the state equation \eqref{Eq:state_weak_form} reads: given $q\in Q_{ad}$, find $w_k(q) \in X_k^0$ such that
\begin{align}\label{Eq:semi_disc_state}
	B(w_k(q),v_k)&=(f,v_k)_I+(u_0,v^{+}_{k,0})-B(q,v_k)\quad \forall v_k \in X_k^0,
\end{align} 
and set $u_k(q) =w_k(q)+q \in X_k^0+Q$ to be the semi-discrete solution of \eqref{Eq:state_weak_form}. 
\begin{remark}
	It is clear that the exact solution $w(q)\in X$ of \eqref{Eq:state_weak_form} satisfies
	\begin{align}
		B(w(q),v_k)&=(f,v_k)_I+(u_0,v^{+}_{k,0})-B(q,v_k)\quad \forall v_k \in X_k^0.
	\end{align} 
	This leads to the Galerkin orthogonality $B(w(q)-w_k(q),v_k)=0\;\text{for all}\; v_k  \in X_k^0$ and hence
	$B(u(q)-u_k(q),v_k)=0\;\text{for all}\; v_k \in X_k^0$.
\end{remark}
A use of integration by parts in time in the bilinear form $B(.,.)$ defined in \eqref{defn:B} yields the equivalent form
\begin{align}\label{Eq:semi_disc_state_1}
	B(w_k,v_k)=-\sum_{m=1}^{M}(w_k,\partial_tv_k)_{I_m}+(\nabla w_k,\nabla v_k)_{I}-\sum_{m=1}^{M-1} \big(w_{k,m}^{-}\jump{v_k}_m\big)+\big(w^{-}_{k,M},v^{-}_{k,M}\big).
\end{align}

\noindent{\bf Discretization in space}.
Let $\cT_h$ be a shape-regular triangulation of $\Omega$. Let $h_K$ denote the diameter of the triangle $K\in\cT_h$ and define the space discretization parameter $h=\max_{K\in\cT_h} h_K$. We consider the conforming finite element space:
\begin{align*}
	V_h:=\{v_h\in H_0^1(\Omega)\;|\;v_h|_{K}\in \mathcal{P}_1(K)\; \text{for} \;K \in \mathcal{T}_h \}.
\end{align*}
Moreover, we consider the fully discrete space-time finite element space
\begin{align*}
	X^{0,1}_{k,h}=\{v_{kh}\in L^2(I,V_h)\;|\;v_{kh}|_{I_m}\in \mathcal{P}_0(I_m,V_h)\}\subseteq X_k^0.
\end{align*}

The fully discrete (space-time discretized) state equation for given control $q\in Q_{ad}$ has the following form: Find $w_{kh}(q)\in X_{k,h}^{0,1} $ such that
\begin{align}\label{Eq:fully_disc_state}
	B(w_{kh}(q),v_{kh})&=(f,v_{kh})_I+(u_0,v^{+}_{kh,0})-B(q,v_{kh})\quad \forall v_{kh} \in X^{0,1}_{k,h},
\end{align}
and set 
$u_{kh}(q):=w_{kh}(q)+q$.
Furthermore, the fully discrete adjoint state equation for given control $q\in Q_{ad}$ has the following form: Find $\phi_{kh}(\bar{q})\in X^{0,1}_{k,h}$ such that
\begin{align}\label{Eq:fully_disc_adj_state}
	B(v_{kh},\phi_{kh}(q),)=\big(u_{kh}(q)-u_d,v_{kh}\big)_I\quad \forall v_{kh} \in X^{0,1}_{k,h}.
\end{align}
We state below the stability result for the fully discretized solutions of the state and adjoint state equations (see \cite[Theorem 4.6 and Corollary 4.7]{Meid_Vexler_08_I}) as:
\begin{lemma}\label{stability_estimate_kh}
	For $q \in Q$, let the solutions $w_{kh}(q)$ and $\phi_{kh}(q)$ be given by the discrete state equation \eqref{Eq:fully_disc_state} and adjoint equation \eqref{Eq:fully_disc_adj_state}, respectively. Then it holds that 
	\begin{align}
		\norm{w_{kh}(q)}_{I}+\norm{\nabla w_{kh}(q)}_{I}&\leq C (\norm{f}_{I}+|q|_{1,\Omega\times I}+\norm{\Pi_h u_{0}}_{0,\Omega}+\norm{\nabla\Pi_h u_{0}}_{0,\Omega}),\label{stability_estimate_state}\\
		\norm{\phi_{kh}(q)}_{I}+\norm{\nabla \phi_{kh}(q)}_{I}&\leq C \norm{u_{kh}(q)-u_d}_{I},\label{stability_estimate_costate}
	\end{align}
	where $\Pi_h:H^1_0(\Omega)\rightarrow V_h$ denotes the spatial $L^2-$projection.
\end{lemma}

\subsection{Error estimates of uncontrolled state and adjoint state variables} This section is devoted to the derivation of \textit{a priori} error estimations for the discrete solutions of the uncontrolled state and adjoint state equation. We introduce some auxiliary equations which are used to simplify our error analysis. For given control $\bar{q}\in Q,$ let $w_{kh}(\bar{q})\in X^{0,1}_{k,h}$ be the fully discrete solution of the following auxiliary state equation:
\begin{align}\label{Aux_state}
	B(w_{kh}(\bar{q}),v_{kh})=(f,v_{kh})_I+(u_0,v^{+}_{kh,0})-B(\bar{q},v_{kh})\quad \forall v_{kh} \in X^{0,1}_{k,h},
\end{align}
and set $u_{kh}(\bar{q}):=w_{kh}(\bar{q})+\bar{q}$. Furthermore, for given $\bar{u}(\bar{q})\in L^2(I,L^2(\Omega)),$ let $\phi_{kh}(\bar{q})\in X^{0,1}_{k,h}$ be the fully discrete solution of the following auxiliary adjoint state equation:
\begin{align}\label{Aux_adj_state}
	B(v_{kh},\phi_{kh}(\bar{q}),)=\big(\bar{u}(\bar{q})-u_d,v_{kh}\big)_I\quad \forall v_{kh} \in X^{0,1}_{k,h}.
\end{align}

\subsubsection{Error estimates of uncontrolled state}
For given fixed control $\bar{q}$ we derive the error between $\bar{u}(\bar{q})$ and $u_{kh}(\bar{q}),$ where $\bar{u}(\bar{q})$ be the solution of the state equation \eqref{eq:state_0a:para}-\eqref{eq:state_a:para} and $u_{kh}(\bar{q})$ be the solution of the fully discrete auxiliary state equation \eqref{Aux_state}. 
Let $u_k(\bar{q})$ be the solution of semi-discrete state equation \eqref{Eq:semi_disc_state} for the given control $\bar{q}$. Since, $\bar{u}(\bar{q})$ and $u_{kh}(\bar{q})$ are state and auxiliary discrete state solution, by the splittings we have $\bar{u}(\bar{q})=\bar{w}(\bar{q})+\bar{q}$ and $u_{kh}(\bar{q})=w_{kh}(\bar{q})+\bar{q}$. Now, we can write the total error $$\bar{u}(\bar{q})-u_{kh}(\bar{q})=\bar{w}(\bar{q})-w_{kh}(\bar{q}).$$
The influences of the
space and time discretization is separated by the temporal part $e_k:=\bar{w}(\bar{q})-w_k(\bar{q})$ and the spatial part $e_h:=w_k(\bar{q})-w_{kh}(\bar{q})$, i.e., $$\bar{w}(\bar{q})-w_{kh}(\bar{q})=e_k+e_h.$$ Our aim is to find the following energy error estimate of the state:
\begin{align}\label{Estm:18}
	\norm{\nabla(\bar{u}(\bar{q})-u_{kh}(\bar{q}))}_I&=\norm{\nabla(\bar{w}(\bar{q})-w_{kh}(\bar{q}))}_I\leq \norm{\nabla e_k}_{I}+\norm{\nabla e_h}_{I}.
\end{align}

\begin{theorem}\label{temporal_error}
	There holds,
	\begin{align*}
		\norm{\nabla e_k}_I\leq C k \norm{\bar{w}(\bar{q})}_{H^1(I;H_0^1(\Omega))},
	\end{align*}
	where $C$ is a positive constant independent of the time step $k$.
\end{theorem}
The proof follows from the following lemmas.
Define the semi-discrete projection \cite{Meid_Vexler_08_I}
\begin{align}
	I_k:C(I,H^1_0(\Omega))\rightarrow X^0_k	
\end{align}
by $I_ku|_{I_m}\in P_0(I_m,H^1_0(\Omega))$ and $I_ku(t_m^{-})=u(t_m^{-})$ for $m=1,2,3,...,M$.
Introducing the projection we get, $$e_k=\bar{w}(\bar{q})-w_k(\bar{q})=(\bar{w}(\bar{q})-I_k\bar{w}(\bar{q}))+(I_k\bar{w}(\bar{q})-w_k(\bar{q}))=\eta_k+\zeta_k,$$
where $\eta_k:=\bar{w}(\bar{q})-I_k\bar{w}(\bar{q})$ and $\zeta_k:=I_k\bar{w}(\bar{q})-w_k(\bar{q})$.
Now, we need to prove the following results: 
\begin{lemma}\label{B_estm}
	The projection error $\eta_k=\bar{w}(\bar{q})-I_k\bar{w}(\bar{q}),$ satisfies
	\begin{align*}
		B(\eta_k,\psi)=(\nabla\eta_k,\nabla \psi)_{I}\quad \forall \psi \in X^0_k.
	\end{align*}	
\end{lemma}
\begin{proof}
	By means of \eqref{Eq:semi_disc_state_1}, we have 
	\begin{align}
		B(\eta_k,\psi)=-\sum_{m=1}^{M}(\eta_k,\partial_t\psi)_{I_m}+(\nabla\eta_k,\nabla \psi)_{I}-\sum_{m=1}^{M-1} \big(\eta_{k,m}^{-}\jump{\psi}_m\big)+\big(\eta^{-}_{k,M},\psi^{-}_{k,M}\big).
	\end{align}
	The term $\eta_{k,m}^{-}$ and $\eta_{k,M}^{-}$ vanishes, because of the definition of interpolation $I_k$. Since $\psi$ lies in the semi-discrete space $X^0_k$, the first term vanishes. The only remaining term is the second one. This completes the proof.
\end{proof} 
\begin{lemma}\label{estm_grad} The temporal error $e_k$  is estimated by the
	projection error $\eta_k$ as
	$$\norm{\nabla e_k}_{I}\leq C\norm{\nabla \eta_k}_{I}.$$
\end{lemma}
\begin{proof}
	Let $\tilde{w}_k\in X^0_k$ be the solution of 
	\begin{align}\label{AuX:Equn_1}
		B(v,\tilde{w}_k)=(\nabla v,\nabla e_k)_{I}\quad \forall v \in X^0_k,
	\end{align}
	with the stability estimate $\norm{\nabla \tilde{w}_k }_{I}\leq C \norm{\nabla e_k}_{I}$.
	By the Galerkin orthogonality, we have
	\begin{align}
		B(\bar{w}(\bar{q})-w_k(\bar{q}),z_k)&=0\quad z_k\in  X^0_k,\nonumber\\
		B(\zeta_k+\eta_k,z_k)&=0\quad z_k\in  X^0_k.\label{Galerkin_ortho}
	\end{align} 
	Using above the stability estimate, Galerkin orthogonality \eqref{Galerkin_ortho}, Lemma \ref{B_estm} and \eqref{AuX:Equn_1}, we obtain the following estimate:
	\begin{align}
		\norm{\nabla e_k}_I^2=&(\nabla e_k,\nabla e_k)_I=(\nabla \eta_k,\nabla e_k)_I+(\nabla \zeta_k,\nabla e_k)_I\nonumber\\
		&=(\nabla \eta_k,\nabla e_k)_I+B(\zeta_k,\tilde{w}_k)=(\nabla \eta_k,\nabla e_k)_I-B(\eta_k,\tilde{w}_k)\nonumber\\
		&=(\nabla \eta_k,\nabla e_k)_{I}-(\nabla \eta_k,\nabla \tilde{w}_k)_I\nonumber\\
		&\leq \norm{\nabla \eta_k}_{I}\norm{\nabla e_k}_I+C\norm{\nabla \eta_k}_{I}\norm{\nabla e_k}_I.
	\end{align}
	Hence, $\norm{\nabla e_k}_I\leq C \norm{\nabla \eta_k}_{I}$. 
\end{proof}
The next lemma on interpolation estimation follows from \cite{Book:Vidar_Thomee}.
\begin{lemma}\label{proj_estm}
	The projection error $\eta_k=\bar{w}(\bar{q})-I_k\bar{w}(\bar{q})$ has the following estimate:
	\begin{align*}
		\norm{\nabla \eta_k}_{I}=\norm{\nabla(\bar{w}(\bar{q})-I_k\bar{w}(\bar{q}))}_I\leq k \norm{\bar{w}(\bar{q})}_{H^1(I;H^1(\Omega))}.
	\end{align*} 
\end{lemma}
\begin{proof}[Proof of Theorem \ref{temporal_error}:]
	Using Lemma \ref{estm_grad} and Lemma \ref{proj_estm}, the proof follows.
\end{proof}

\begin{theorem}\label{spatial_error}
	Let $w_k(\bar{q})$  be the semidiscretized
	solution  of \eqref{Eq:semi_disc_state} and $w_{kh}(\bar{q})$ be the fully discretized solution of \eqref{Eq:fully_disc_state}. Then the error $e_h=w_k(\bar{q})-w_{kh}(\bar{q})$ has the estimate
	\begin{align*}
		\norm{\nabla e_h}_I\leq C h^{\tau-1} \norm{w_k}_{L^2(I;H^{\tau}(\Omega))},
	\end{align*}	
	where the constant $C$ is independent of the mesh size $h$ and the size of the time steps $k$.
\end{theorem}
The proof is divided into several steps which
are collected in the following lemmas. Define the projection $$\pi_h:X^0_k\rightarrow X^{0,1}_{k,h}$$ by means of the spatial $L^2$-projection $\Pi_h:H^1_0(\Omega)\rightarrow V_h$ point-wise in time as $(\pi_h w_k)(t)=\Pi_hw_k(t)$.
Introducing the projection $\pi_h$ we get, $$e_h=(w_k(\bar{q})-\pi_hw_k(\bar{q}))+(\pi_hw_k(\bar{q})-w_{kh}(\bar{q}))=\eta_h+\zeta_h,$$
where $\eta_h:=w_k(\bar{q})-\pi_hw_k(\bar{q})$ and $\zeta_h:=\pi_hw_k(\bar{q})-w_{kh}(\bar{q})$. 
\begin{lemma}
	The projection error $\eta_h$ satisfies the following relation
	\begin{align*}
		B(\eta_h,\psi)=(\nabla \eta_h,\nabla \psi)_{I} \quad\forall \psi \in X^{0,1}_{k,h}.
	\end{align*}
\end{lemma}
\begin{proof}
	The proof follows by similar arguments of Lemma \ref{B_estm}. 
\end{proof}
\begin{lemma}\label{stability_estm} The following boundedness property of the error $\zeta_h$ holds
	\begin{align*}
		\norm{\nabla \zeta_h}_{I}\leq C\norm{\nabla \eta_h}_{I}.
	\end{align*}
\end{lemma}
\begin{proof}
	For $v\in X^0_k$, the definition of $B$ in \eqref{defn:B} reads
	\begin{align}\label{defn:B1}
		B(v,v)=\sum_{m=1}^{M}\big(\partial_t v,v\big)_{I_m}+\big(\nabla v , \nabla v\big)_{I}+\sum_{m=1}^{M-1} \big(\jump{v}_m,v^{+}_{m}\big)+\big(v^{+}_{0},v^{+}_{0}\big),
	\end{align}
	and in \eqref{Eq:semi_disc_state_1} as
	\begin{align}\label{defn:B2}
		B(v,v)=-\sum_{m=1}^{M}\big(v,\partial_t v\big)_{I_m}+\big(\nabla v , \nabla v\big)_{I}-\sum_{m=1}^{M-1} \big(v^{-}_{m},\jump{v}_m\big)+\big(v^{-}_{M},v^{-}_{M}\big).
	\end{align}	
	Adding the above two equations \eqref{defn:B1} and \eqref{defn:B2}, we get
	\begin{align}\label{Estm:39}
		B(v,v)\geq \big(\nabla v , \nabla v\big)_{I} \;\text{for all}\quad v\in X^0_k.
	\end{align}
	Choosing, $v=\zeta_h$ in \eqref{Estm:39} and utilizing the Galerkin orthogonality of the space discretization, we obtain
	$$\norm{\nabla \zeta_h}_{I}^2\leq B(\zeta_h,\zeta_h)=-B(\eta_h,\zeta_h)=-(\nabla \eta_h,\nabla \zeta_h)\leq \norm{\nabla \eta_h}_I\norm{\nabla\zeta_h}_I.$$
	This establishes the desired result $\norm{\nabla \zeta_h}_{I}\leq \norm{\nabla \eta_h}_{I}$.
\end{proof}
We state below the well-known estimate for the spatial projection  $\pi_h$ (see \cite[subsection 5.2]{Meid_Vexler_08_I}) as:
\begin{lemma}\label{Interpolation_estm}
	The projection error $\eta_h=w_k(\bar{q})-\pi_hw_k(\bar{q})$ has the following estimate: 
	\begin{align}
		\norm{\nabla(w_k(\bar{q})-\pi_hw_k(\bar{q}))}_I\leq C h^{\tau-1} \norm{w_k}_{L^2(I;H^{\tau}(\Omega))}.
	\end{align}
\end{lemma}
\noindent
\begin{proof}[Proof of Theorem \ref{spatial_error}:]
	The splitting $e_h=\eta_h+\zeta_h$ yields $\norm{\nabla e_h}_I\leq \norm{\nabla \eta_h}_I+\norm{\nabla \zeta_h}_I$. Using Lemma \ref{stability_estm}, we obtain  
	$\norm{\nabla e_h}_I\leq C \norm{\nabla \eta_h}_I$. Finally, using Lemma \ref{Interpolation_estm}, we get the desired result
	$$\norm{\nabla e_h}_I\leq C h^{\tau-1} \norm{w_k}_{L^2(I;H^{\tau}(\Omega))}.$$	
\end{proof} 
The next theorem follows from Theorem \ref{temporal_error} and \ref{spatial_error}. 
\begin{theorem}\label{Error_estm_state_fixed_q} 
	The error estimation for the uncontrolled state variable is given by
	\begin{align}
		\norm{\nabla(\bar{u}(\bar{q})-u_{kh}(\bar{q}))}_I=\norm{\nabla(\bar{w}(\bar{q})-w_{kh}(\bar{q}))}_I\leq  C \big( k +h^{\tau-1} \big).\nonumber
	\end{align}
\end{theorem}



Now, we derive the error estimate of the uncontrolled adjoin state i.e., for a given fixed control $\bar{q}$ we derive the error between $\bar{\phi}(\bar{q})$ and $\phi_{kh}(\bar{q})$, where $\bar{\phi}(\bar{q})$ be the solution of \eqref{eq:adjstate_a:para} and  $\phi_{kh}(\bar{q})$ be the solution of \eqref{Eq:fully_disc_adj_state} for $q=\bar{q}$.
\begin{theorem}\label{estm_adj_fixed_q} There holds,
	\begin{align*}
		\norm{\nabla(\bar{\phi}(\bar{q})-\phi_{kh}(\bar{q}))}_I&\leq C( k+ h^{\tau-1}).	
	\end{align*}
\end{theorem}
\begin{proof}
	Let $\tilde{\phi}_k(\bar{q})\in X^0_k$ be the solution of the auxiliary problem:
	\begin{align*}
		B(v_k,\tilde{\phi}_k(\bar{q}))=(v_k,\bar{u}(\bar{q})-u_d)_I \quad \forall v_k\in X^0_k.
	\end{align*} 
	Introducing $\tilde{\phi}_k(\bar{q})$ and applying the triangle inequality, we arrive at
	\begin{align}\label{Estm:25}
		\norm{\nabla(\bar{\phi}(\bar{q})-\phi_{kh}(\bar{q}))}_I&\leq \norm{\nabla(\bar{\phi}(\bar{q})-\tilde{\phi}_k(\bar{q}))}_I+\norm{\nabla(\tilde{\phi}_k(\bar{q})-\phi_{kh}(\bar{q}))}_I.
	\end{align}
	For the first term of \eqref{Estm:25}, we apply similar arguments of the proof of Theorem \ref{temporal_error} to obtain
	\begin{align}\label{Estm:23}
		\norm{\nabla(\bar{\phi}(\bar{q})-\tilde{\phi}_k(\bar{q}))}_I\leq C k.
	\end{align}
	Using the stability result \eqref{stability_estimate_costate}, Poincar\'e inequality, and Theorem \ref{Error_estm_state_fixed_q} we obtain 
	\begin{align}\label{Estm:24}
		\norm{\nabla(\tilde{\phi}_k(\bar{q})-\phi_{kh}(\bar{q}))}_I\leq C \norm{\nabla(\bar{w}(\bar{q})-w_{kh}(\bar{q}))}_I \leq C h^{\tau-1}.	
	\end{align}
	Putting the estimates \eqref{Estm:23} and \eqref{Estm:24} in \eqref{Estm:25}, we get the required result	
	
\end{proof}

\subsection{Discretization of the control variable}
In this subsection, we describe the discretization of the control variable. Let 
$Q_{\sigma}$ be an admissible discrete subspace of the control space $Q$. 
One can consider different mesh for control variable than the mesh corresponding to state and adjoint variables, see \cite{Vexler_et.al:2007}. However,
for simplicity of notation we will use the same time-partitioning \eqref{Time_partition} and the same spatial
mesh $\cT_h$ defined in subsection~\ref{subsec:Discretization}.  

Let $\cT_{\sigma}$ be a regular mesh of the domain $\Omega \times (0,T)$, consists of prism see
\cite{Ern:2004:FEMbook,Ciarlet:1978:FEM}. A typical prism is denoted by $K \times I_m,$ where $K\in \cT_h$ and $I_m$ is a subinterval of $I$ such that $c_1 h_K\leq |I_m|\leq c_2 h_K$ for some fixed positive constants $c_1$ and $c_2$. Using a spatial mesh $\cT_h$ we consider the following finite element space:
\begin{align*}
	Q_h=\{q_h\in C(\bar{\Omega})\;|\;q_h|_{K}\in \mathcal{P}_1(K)\; \text{for} \;K \in \mathcal{T}_h\}.
\end{align*}
Define a conforming $\mathcal{P}_1$-subspace
$Q_{\sigma}$ of $Q$ by,
\begin{align}\label{discrete_control_space}
	Q_{\sigma}=\{q_{\sigma}\in C(\bar{\Omega}\times\bar{I})\;|\;q_{\sigma}|_{K \times I_m}\in \mathcal{P}_1(K)\otimes \mathcal{P}_1(I_m),\; q_{\sigma}(x,0)=q_{\sigma}(x,T)=0,\; x\in \Omega\}.
\end{align}
Another equivalent representation of $Q_{\sigma}$ is the following: 
\begin{align}
	Q_{\sigma}=\{q_{\sigma}\in C(\bar{\Omega}\times\bar{I})\,|\,q_{\sigma}|_{I_m}\in \mathcal{P}_1(I_m,Q_h),\; q_{\sigma}(x,0)=q_{\sigma}(x,T)=0,\; x\in \Omega\}.
\end{align}
The discrete admissible control set,
\begin{align}
	Q^{\sigma}_{ad}=\{q_{\sigma}\in Q_{\sigma}\,|\, q_a\leq q_{\sigma}(p_{\sigma
	})\leq q_b \quad \text{for all} \; \text{nodes}\; p_{\sigma
	}\; \text{on}\; \partial \Omega \times (0,T)\}.
\end{align}
The fully discrete optimal control problem is given as follows:
\begin{align}\label{Min:non_redu_control_disc}
	\text{Minimize}\; J(u_{kh}(q_{\sigma}),q_{\sigma})\; \text{subject to \eqref{Eq:fully_disc_state} and}\; (u_{kh}(q_{\sigma}),q_{\sigma})\in  (X^{0,1}_{k,h}+Q_{\sigma})\times Q^{\sigma}_{ad}.
\end{align}
Here, the parameter $\sigma$ collects the discretization parameters $k$ and $h$ i.e., $\sigma=\sigma(k,h)$. 
The standard theory of optimal control problem \cite{Raymond:Notes,trolzstch:2005:Book} can be employed to deduce the existence and uniqueness of the solution of the following discrete optimality system:
\begin{proposition}[Fully discrete optimality system]
	There exists a unique solution $(\bar{u}_{kh}(\bar{q}_{\sigma}),\bar{q}_{\sigma})\in (X^{0,1}_{k,h}+Q_{\sigma}) \times Q^{\sigma}_{ad}$ for the Dirichlet control problem  $\eqref{Min:non_redu_control_disc}$. Further, there exists a unique adjoint state $\bar{\phi}_{kh}(\bar{q}_{\sigma})\in X_{k,h}^{0,1}$ satisfies the following:
	\begin{subequations}
		\begin{align}
			\bar{u}_{kh}(\bar{q}_{\sigma})=&\bar{w}_{kh}(\bar{q}_{\sigma})+\bar{q}_{\sigma},\quad \bar{w}_{kh}(\bar{q}_{\sigma}) \in X_{k,h}^{0,1}.\\
			B(\bar{w}_{kh}(\bar{q}_{\sigma}),v_{kh})&=(f,v_{kh})_I+(u_0,v^{+}_{kh,0})-B(\bar{q}_{\sigma},v_{kh})\quad \forall v_{kh} \in X^{0,1}_{k,h}.\label{full_disc_state}\\
			B(v_{kh},\bar{\phi}_{kh}(\bar{q}_{\sigma}),)&=\big(\bar{u}_{kh}(\bar{q}_{\sigma})-u_d,v_{kh}\big)_I\quad \forall v_{kh} \in X^{0,1}_{k,h}.\label{full_disc_adj_state}\\
			\lambda(\partial_t \bar{q}_{\sigma},\partial_t (p_{\sigma}-\bar{q}_{\sigma}))_{I}+&\lambda (\nabla \bar{q}_{\sigma},\nabla(p_{\sigma}-\bar{q}_{\sigma}))_I\geq  (\partial_t(p_{\sigma}-\bar{q}_{\sigma}), \bar{\phi}_{kh}(\bar{q}_{\sigma}))_I\nonumber\\&+(\nabla \bar{\phi}_{kh}(\bar{q}_{\sigma}),\nabla(p_{\sigma}-\bar{q}_{\sigma}))_I-(\bar{u}_{kh}(\bar{q}_{\sigma})-u_d,p_{\sigma}-\bar{q}_{\sigma})_I, \label{Discrete_VI}
		\end{align}	
	\end{subequations}
for all $p_{\sigma} \in Q^{\sigma}_{ad}.$
\end{proposition}

\noindent
The next lemma is used for the derivation of error estimate below.
\begin{lemma}\label{aux_lemma} There holds,
	\begin{align}
		(\partial_t(\bar{q}_{\sigma}-\bar{q}),\bar{\phi}_{kh}(\bar{q}_{\sigma})-\bar{\phi}_{kh}(\bar{q}))_I+&(\nabla(\bar{q}_{\sigma}-\bar{q}),\nabla(\bar{\phi}_{kh}(\bar{q}_{\sigma})-\bar{\phi}_{kh}(\bar{q})))_I\\&=\big(\bar{u}_{kh}(\bar{q}_{\sigma})-\bar{u}(\bar{q}),w_{kh}(\bar{q})-\bar{w}_{kh}(\bar{q}_{\sigma})\big)_I.
	\end{align}
\end{lemma}
\begin{proof}
	A subtraction of \eqref{full_disc_state} from \eqref{Aux_state} yields, 
	\begin{align}\label{Aux_inter}
		B(w_{kh}(\bar{q})-\bar{w}_{kh}(\bar{q}_{\sigma}),v_{kh})=B(\bar{q}_{\sigma}-\bar{q},v_{kh})\quad \forall v_{kh} \in X^{0,1}_{k,h}.
	\end{align}
	Also, subtracting, \eqref{Aux_adj_state} from  \eqref{full_disc_adj_state} we get, 
	\begin{align}\label{Aux_inter_1}
		B(v_{kh},\bar{\phi}_{kh}(\bar{q}_{\sigma})-\bar{\phi}_{kh}(\bar{q}),)=\big(\bar{u}_{kh}(\bar{q}_{\sigma})-\bar{u}(\bar{q}),v_{kh}\big)_I\quad \forall v_{kh} \in X^{0,1}_{k,h}.
	\end{align}
	Choose, $v_{kh}=\bar{\phi}_{kh}(\bar{q}_{\sigma})-\bar{\phi}_{kh}(\bar{q})$ in \eqref{Aux_inter} we get,
	\begin{align}\label{estm1}
		B(w_{kh}(\bar{q})-\bar{w}_{kh}(\bar{q}_{\sigma}),\bar{\phi}_{kh}(\bar{q}_{\sigma})-\bar{\phi}_{kh}(\bar{q}))=B(\bar{q}_{\sigma}-\bar{q},\bar{\phi}_{kh}(\bar{q}_{\sigma})-\bar{\phi}_{kh}(\bar{q})).
	\end{align}
	Choose, $v_{kh}=w_{kh}(\bar{q})-\bar{w}_{kh}(\bar{q}_{\sigma})$ in \eqref{Aux_inter_1}, we obtain
	\begin{align}\label{estm3}
		B(w_{kh}(\bar{q})-\bar{w}_{kh}(\bar{q}_{\sigma}),\bar{\phi}_{kh}(\bar{q}_{\sigma})-\bar{\phi}_{kh}(\bar{q}),)=\big(\bar{u}_{kh}(\bar{q}_{\sigma})-\bar{u}(\bar{q}),w_{kh}(\bar{q})-\bar{w}_{kh}(\bar{q}_{\sigma})\big)_I. 
	\end{align}
	Now, equating \eqref{estm1} and \eqref{estm3}, we have
	\begin{align}
		B(\bar{q}_{\sigma}-\bar{q},\bar{\phi}_{kh}(\bar{q}_{\sigma})-\bar{\phi}_{kh}(\bar{q}))=\big(\bar{u}_{kh}(\bar{q}_{\sigma})-\bar{u}(\bar{q}),w_{kh}(\bar{q})-\bar{w}_{kh}(\bar{q}_{\sigma})\big)_I.
	\end{align} 
	Computing $B(\bar{q}_{\sigma}-\bar{q},\bar{\phi}_{kh}(\bar{q}_{\sigma})-\bar{\phi}_{kh}(\bar{q})),$ we have
	\begin{align*}
		(\partial_t(\bar{q}_{\sigma}-\bar{q}),\bar{\phi}_{kh}(\bar{q}_{\sigma})-\bar{\phi}_{kh}(\bar{q}))_I+&(\nabla(\bar{q}_{\sigma}-\bar{q}),\nabla(\bar{\phi}_{kh}(\bar{q}_{\sigma})-\bar{\phi}_{kh}(\bar{q})))_I\\&=\big(\bar{u}_{kh}(\bar{q}_{\sigma})-\bar{u}(\bar{q}),w_{kh}(\bar{q})-\bar{w}_{kh}(\bar{q}_{\sigma})\big)_I.
	\end{align*} 	
\end{proof}

\begin{theorem}[Error estimation for control]\label{Estm:Control} There holds,
	\begin{align}\label{Hilld_term}
		\lambda |\bar{q}-\bar{q}_{\sigma}|^2_{1,\Omega\times I}+&\norm{\bar{u}(\bar{q})-\bar{u}_{kh}(\bar{q}_{\sigma})}^2_I
		\leq 
		\big[\lambda(\partial_t \bar{q},\partial_t (p_{\sigma}-\bar{q}))_{I}+\lambda (\nabla \bar{q},\nabla (p_{\sigma}-\bar{q}))_I\nonumber\\& 
		-(\partial_t(p_{\sigma}-\bar{q}),\bar{\phi})_I-(\nabla(p_{\sigma}-\bar{q}),\nabla\bar{\phi}(\bar{q}))_I+(\bar{u}(\bar{q})-u_d,p_{\sigma}-\bar{q})_I\big]\nonumber\\&+\norm{\bar{q}-p_{\sigma}}^2_{1,\Omega\times I}
		+\norm{\nabla(\bar{\phi}(\bar{q})-\phi_{kh}(\bar{q}))}^2_I+\norm{\bar{w}(\bar{q})-w_{kh}(\bar{q})}^2_I
	\end{align}
	for all $p_{\sigma}\in Q^{\sigma}_{ad}$.
\end{theorem}
\begin{proof}
	Choose, $p=\bar{q}_{\sigma}$ in \eqref{conti_VI:para}, we get 
	\begin{align}\label{Auxi_Cont_VI}
		\lambda(\partial_t \bar{q},\partial_t (\bar{q}_{\sigma}-\bar{q}))_{I}+&\lambda (\nabla \bar{q}, \nabla (\bar{q}_{\sigma}-\bar{q}))_I\geq -(\bar{\phi}(\bar{q}),\partial_t(\bar{q}_{\sigma}-\bar{q}))_I\nonumber\\&+(\nabla \bar{\phi}(\bar{q}), \nabla (\bar{q}_{\sigma}-\bar{q}))_I-(\bar{u}(\bar{q})-u_d,\bar{q}_{\sigma}-\bar{q})_I.	
	\end{align}
	Rearranging the terms for the discrete variational inequality \eqref{Discrete_VI}, we get
	\begin{align}\label{Auxi_disc_VI}
		\lambda(\partial_t \bar{q}_{\sigma},\partial_t (\bar{q}-\bar{q}_{\sigma}))_{I}+&\lambda (\nabla \bar{q}_{\sigma},\nabla (\bar{q}-\bar{q}_{\sigma}))_I\geq -\lambda(\partial_t \bar{q}_{\sigma},\partial_t (p_{\sigma}-\bar{q}))_{I}\nonumber\\&-\lambda (\nabla \bar{q}_{\sigma},\nabla (p_{\sigma}-\bar{q}))_I  +(\partial_t(p_{\sigma}-\bar{q}),\bar{\phi}_{kh}(\bar{q}_{\sigma}))_I\nonumber\\&+(\partial_t(\bar{q}-\bar{q}_{\sigma}),\bar{\phi}_{kh}(\bar{q}_{\sigma}))_I+(\nabla(p_{\sigma}-\bar{q}),\nabla\bar{\phi}_{kh}(\bar{q}_{\sigma}))_I\nonumber\\&+(\nabla(\bar{q}-\bar{q}_{\sigma}),\nabla\bar{\phi}_{kh}(\bar{q}_{\sigma}))_I-(\bar{u}_{kh}(\bar{q}_{\sigma})-u_d,p_{\sigma}-\bar{q})_I\nonumber\\&-(\bar{u}_{kh}(\bar{q}_{\sigma})-u_d,\bar{q}-\bar{q}_{\sigma})_I	
	\end{align}
	for all $p_{\sigma}\in Q^{\sigma}_{ad}$.
	Adding \eqref{Auxi_Cont_VI} and \eqref{Auxi_disc_VI}, we get 
	\begin{align}
		-\lambda \norm{\partial_t  (\bar{q}-\bar{q}_{\sigma})}^2_{I}-&\lambda \norm{\nabla (\bar{q}-\bar{q}_{\sigma})}^2_I
		\geq 
		\lambda(\partial_t \bar{q}_{\sigma},\partial_t (\bar{q}-p_{\sigma}))_{I}+\lambda (\nabla \bar{q}_{\sigma},\nabla (\bar{q}-p_{\sigma}))_I\nonumber\\&  +(\partial_t(p_{\sigma}-\bar{q}),\bar{\phi}_{kh}(\bar{q}_{\sigma}))_I+(\partial_t(\bar{q}-\bar{q}_{\sigma}),\bar{\phi}_{kh}(\bar{q}_{\sigma}))_I\nonumber\\&
		+(\nabla(p_{\sigma}-\bar{q}),\nabla\bar{\phi}_{kh}(\bar{q}_{\sigma}))_I+(\nabla(\bar{q}-\bar{q}_{\sigma}),\nabla\bar{\phi}_{kh}(\bar{q}_{\sigma}))_I\nonumber\\&
		-(\bar{u}_{kh}(\bar{q}_{\sigma})-u_d,p_{\sigma}-\bar{q})_I-(\bar{u}_{kh}(\bar{q}_{\sigma})-u_d,\bar{q}-\bar{q}_{\sigma})_I\nonumber\\
		&\geq  
		\big[\lambda(\partial_t \bar{q},\partial_t (\bar{q}-p_{\sigma}))_{I}+\lambda (\nabla \bar{q},\nabla (\bar{q}-p_{\sigma}))_I\nonumber\\& 
		+(\partial_t(p_{\sigma}-\bar{q}),\bar{\phi}(\bar{q}))_I+(\nabla(p_{\sigma}-\bar{q}),\nabla\bar{\phi}(\bar{q}))_I\nonumber\\&-(\bar{u}(\bar{q})-u_d,p_{\sigma}-\bar{q})_I\big]+\lambda(\partial_t( \bar{q}_{\sigma}-\bar{q}),\partial_t (\bar{q}-p_{\sigma}))_{I}\nonumber\\&+\lambda (\nabla (\bar{q}_{\sigma}-\bar{q}),\nabla (\bar{q}-p_{\sigma}))_I
		+(\partial_t(p_{\sigma}-\bar{q}),\bar{\phi}_{kh}(\bar{q}_{\sigma})-\bar{\phi}(\bar{q}))_I\nonumber\\&+(\nabla(p_{\sigma}-\bar{q}),\nabla(\bar{\phi}_{kh}(\bar{q}_{\sigma})-\bar{\phi}(\bar{q})))_I
		-(\bar{u}_{kh}(\bar{q}_{\sigma})-\bar{u}(\bar{q}),p_{\sigma}-\bar{q})_I\nonumber\\&+(\nabla(\bar{q}-\bar{q}_{\sigma}),\nabla(\bar{\phi}_{kh}(\bar{q}_{\sigma})-\bar{\phi}(\bar{q})))_I
		-(\bar{u}_{kh}(\bar{q}_{\sigma})-\bar{u}(\bar{q}),\bar{q}-\bar{q}_{\sigma})_I\nonumber\\&+(\partial_t(\bar{q}-\bar{q}_{\sigma}),\bar{\phi}_{kh}(\bar{q}_{\sigma})-\bar{\phi}(\bar{q}))_I\label{estm5} 	
	\end{align}
	for all $p_{\sigma}\in Q^{\sigma}_{ad}$. Now we need to do some manipulation on the last three terms in \eqref{estm5}. Denote  
	\begin{align}\label{Estm:9}
		E=(\nabla(\bar{q}-\bar{q}_{\sigma}),\nabla(\bar{\phi}_{kh}(\bar{q}_{\sigma})-\bar{\phi}(\bar{q})))_I-(\bar{u}_{kh}(\bar{q}_{\sigma})-\bar{u}(\bar{q}),\bar{q}-\bar{q}_{\sigma})_I\nonumber\\+(\partial_t(\bar{q}-\bar{q}_{\sigma}),\bar{\phi}_{kh}(\bar{q}_{\sigma})-\bar{\phi}(\bar{q}))_I.
	\end{align}
	Introducing the auxiliary solution $\phi_{kh}(\bar{q})$ in the first, third term and modifying the second term in \eqref{Estm:9}, we obtain
	\begin{align}
		E=&(\nabla(\bar{q}-\bar{q}_{\sigma}),\nabla(\bar{\phi}_{kh}(\bar{q}_{\sigma})-\phi_{kh}(\bar{q})))_I
		+\norm{\bar{u}_{kh}(\bar{q}_{\sigma})-\bar{u}(\bar{q})}^2_I\nonumber\\&
		-(\bar{u}_{kh}(\bar{q}_{\sigma})-\bar{u}(\bar{q}),\bar{w}_{kh}\bar{q}_{\sigma}-\bar{w}(\bar{q}))_I
		+(\partial_t(\bar{q}-\bar{q}_{\sigma}),\bar{\phi}_{kh}(\bar{q}_{\sigma})-\phi_{kh}(\bar{q}))_I\nonumber\\&
		+(\nabla(\bar{q}-\bar{q}_{\sigma}),\nabla(\phi_{kh}(\bar{q})-\bar{\phi}(\bar{q})))_I\nonumber\\&
		+(\partial_t(\bar{q}-\bar{q}_{\sigma}),\bar{\phi}_{kh}(\bar{q}_{\sigma})-\phi_{kh}(\bar{q}))_I.\label{Estm:10}
	\end{align}
	Using the Lemma \ref{aux_lemma} in \eqref{Estm:10}, we get 
	\begin{align}
		E=&\norm{\bar{u}_{kh}(\bar{q}_{\sigma})-\bar{u}(\bar{q})}^2_I+(\nabla(\bar{q}-\bar{q}_{\sigma}),\nabla(\phi_{kh}(\bar{q})-\bar{\phi}(\bar{q})))_I\nonumber\\&+(\partial_t(\bar{q}-\bar{q}_{\sigma}),\bar{\phi}_{kh}(\bar{q}_{\sigma})-\phi_{kh}(\bar{q}))_I
		-(w_{kh}(\bar{q})-w_{kh}(\bar{q}_{\sigma}),\bar{u}_{kh}(\bar{q}_{\sigma})-\bar{u}(\bar{q}))_I\nonumber\\&
		-(\bar{u}_{kh}(\bar{q}_{\sigma})-\bar{u}(\bar{q}),\bar{w}_{kh}(\bar{q}_{\sigma})-\bar{w}(\bar{q}))_I.\nonumber	
	\end{align}
	Hence,
	\begin{align}\label{estm4}
		E=&\norm{\bar{u}_{kh}(\bar{q}_{\sigma})-\bar{u}(\bar{q})}^2_I+(\nabla(\bar{q}-\bar{q}_{\sigma}),\nabla(\phi_{kh}(\bar{q})-\bar{\phi}(\bar{q})))_I\nonumber\\&+(\partial_t(\bar{q}-\bar{q}_{\sigma}),\bar{\phi}_{kh}(\bar{q}_{\sigma})-\phi_{kh}(\bar{q}))_I
		-(\bar{u}_{kh}(\bar{q}_{\sigma})-\bar{u}(\bar{q}),w_{kh}(\bar{q})-\bar{w}(\bar{q}))_I.	
	\end{align}
	
	\noindent
	Using \eqref{estm4} in \eqref{estm5} and grouping the terms, we get
	\begin{align}
		\lambda \norm{\partial_t  (\bar{q}-\bar{q}_{\sigma})}^2_{I}+&\lambda \norm{\nabla (\bar{q}-\bar{q}_{\sigma})}^2_I+\norm{\bar{u}_{kh}(\bar{q}_{\sigma})-\bar{u}(\bar{q})}^2_I
		\leq 
		\big[\lambda(\partial_t \bar{q},\partial_t (p_{\sigma}-\bar{q}))_{I}\nonumber\\&+\lambda (\nabla \bar{q},\nabla (p_{\sigma}-\bar{q}))_I 
		-(\partial_t(p_{\sigma}-\bar{q}),\bar{\phi}(\bar{q}))_I-(\nabla(p_{\sigma}-\bar{q}),\nabla\bar{\phi}(\bar{q}))_I\nonumber\\&+(\bar{u}(\bar{q})-u_d,p_{\sigma}-\bar{q})_I\big]+\lambda(\partial_t(\bar{q}-\bar{q}_{\sigma}),\partial_t (\bar{q}-p_{\sigma}))_{I}\nonumber\\&+\lambda (\nabla (\bar{q}-\bar{q}_{\sigma}),\nabla (\bar{q}-p_{\sigma}))_I
		+(\partial_t(\bar{q}-p_{\sigma}),\bar{\phi}_{kh}(\bar{q}_{\sigma})-\bar{\phi}(\bar{q}))_I\nonumber\\&+(\nabla(\bar{q}-p_{\sigma}),\nabla(\bar{\phi}_{kh}(\bar{q}_{\sigma})-\bar{\phi}(\bar{q})))_I+(\bar{u}_{kh}(\bar{q}_{\sigma})-\bar{u}(\bar{q}),p_{\sigma}-\bar{q})_I
		\nonumber\\&+(\nabla(\bar{q}-\bar{q}_{\sigma}),\nabla(\bar{\phi}_{kh}(\bar{q})-\bar{\phi}(\bar{q})))_I-(\partial_t(\bar{q}-\bar{q}_{\sigma}),\bar{\phi}_{kh}(\bar{q})-\bar{\phi}(\bar{q}))_I\nonumber\\&-(\bar{u}_{kh}(\bar{q}_{\sigma})-\bar{u}(\bar{q}),w_{kh}(\bar{q})-\bar{q}_{\sigma})_I.\label{Estm:11}	
	\end{align}
	Using the stability estimate of the adjoint state equation, we get
	\begin{align}
		\norm{\nabla(\bar{\phi}_{kh}(\bar{q}_{\sigma})-\bar{\phi}(\bar{q}))}_I&\leq \norm{\nabla(\bar{\phi}_{kh}(\bar{q}_{\sigma})-\phi_{kh}(\bar{q}))}_I+\norm{\nabla(\phi_{kh}(\bar{q})-\bar{\phi}(\bar{q}))}_I\nonumber\\
		&\leq \norm{\bar{u}_{kh}(\bar{q}_{\sigma})-\bar{u}(\bar{q})}_I+\norm{\nabla(\phi_{kh}(\bar{q})-\bar{\phi}(\bar{q}))}_I.\label{Equn:6}
	\end{align}
	Applying the Cauchy--Schwarz inequality and putting \eqref{Equn:6}
	in the above equation \eqref{Estm:11} we obtain the desired estimate.
\end{proof}
Now we derive the convergence rates for the terms on the right-hand side of \eqref{Hilld_term}.  We construct a suitable approximation $p_{\sigma}$ for $\bar{q}$ through some interpolations which are described below.
Let $\mathcal{I}_{\sigma}$ be the Lagrange interpolation operator on the three dimensional prismatic elements. On a prismatic element $K_{\sigma}:=K\times I_m$ with $I_m=(t_{m-1},t_m]$ define the local Lagrange interpolation operator $\mathcal{I}_{K_{\sigma}}$ by the following:
\begin{align}
	\mathcal{I}_{K_{\sigma}}\bar{q}(x,t)=\sum_{i=1}^{3}\Big(\sum_{j=1}^{2}\bar{q}(x_i,t_j)\chi_j(t)\Big)\phi_i(x)
\end{align}
for $(x,t)\in K_{\sigma},$ and $\{\chi_1(t)=(t_m-t)/(t_{m}-t_{m-1}),\,\chi_2(t)=(t-t_{m-1})/(t_{m}-t_{m-1})\}$ temporal basis and $\{\phi_i\}_{i=1}^3$ spatial nodal basis. Let $M_{\sigma}$ be the trace of the discrete control space $Q_{\sigma}$ (see \eqref{discrete_control_space}) on $\Gamma_{C}$, and the discrete extension operator  $R_{\sigma}$ be a map from $M_{\sigma}$ to $Q_{\sigma}$. In \cite{Bernardi:et.al,Scott_et_al}, the discrete extension operator is obtained by combining a standard continuous extension operator with a local regularization operator. Now we define a quasi-interpolation
operator $\mathcal{J}_{\sigma}:W^{1,1}(\Gamma_C)\rightarrow M_{\sigma}$ as follows. Let $v\in W^{1,1}(\Gamma_C)$. For interior
nodes $\mathbf{p}$ in $\Gamma_C$, we choose the Chen--Nochetto operator (see \cite{Chen:Nochetto:2000}) which preserves local affine functions and positivity:
\begin{align*}
	\mathcal{J}_{\sigma}v(\mathbf{p})=\frac{1}{meas(\mathcal{B})}\int_\mathcal{B} v,
\end{align*}
where $\mathcal{B}$ is the largest open ball centered at $\mathbf{p}$ such that it is contained in the union of the elements containing $\mathbf{p}$. For the boundary nodes $\mathbf{p}$ on $\bar{\Gamma}_C\cap \bar{\Gamma}_D,$ we set $\mathcal{J}_{\sigma}v(\mathbf{p})= 0$. For the other boundary nodes $\mathbf{p}$ on $\bar{\Gamma}_C$ we set 
\begin{align*}
	\mathcal{J}_{\sigma}v(\mathbf{p})=\frac{1}{meas(L)}\int_L v,
\end{align*}
where $L$ is a small line segment symmetrically placed
around $\mathbf{p}$, and included in $\bar{\Gamma}_C$. This definition preserves both sign and affine functions. Also, we have the following estimate (see \cite[Corollary 4.2.3]{Book:Ziemer} and \cite{Hild_15_signorini}):
\begin{align}\label{Estm:12}
	\norm{v-\mathcal{J}_{\sigma}v}_{0,K_{\sigma}\cap\Gamma_C}\leq C \norm{\nabla v}_{L^1(K_{\sigma}\cap\Gamma_C)}.
\end{align}
Note that the estimate of the above type \eqref{Estm:12} can not be obtained for the Lagrange interpolation operator $\mathcal{I}_{\sigma}$. Moreover, $\mathcal{J}_{\sigma}$ obeys the same approximation properties as of the Lagrange interpolation.
Now we choose the approximation $p_{\sigma}$ for the control $\bar{q}$ as:
\begin{align}\label{eqn_defn_psig}
	p_{\sigma}=\mathcal{I}_{\sigma}\bar{q}+\mathcal{R}_{\sigma}\big(\mathcal{J}_{\sigma}(\bar{q}|_{\Gamma_C})-\mathcal{I}_{\sigma}(\bar{q}|_{\Gamma_C})\big)\in Q^{\sigma}_{ad}.
\end{align}

To estimate the best approximation term in the bracket of \eqref{Hilld_term}, we introduce the following notations. Let $K_{\sigma}$ be a prism which shares a face with $\bar{\Gamma}_C$. Define $$S_{NC}=\{(x,t)\in K_{\sigma} \cap\Gamma_C:\; q_a<\bar{q}(x,t)<q_b\},$$ and
$$S_{C}=\{(x,t)\in K_{\sigma}\cap\Gamma_C:\; \bar{q}(x,t)=q_a\}\cup\{(x,t)\in K_{\sigma}\cap\Gamma_C:\; \bar{q}(x,t)=q_b\}.$$
The sets $S_{C}$ and $S_{NC}$ are measurable since $q$ is continuous on $\Gamma_C.$ We denote $|S_{C}|$ and $|S_{NC}|$ are their measures. We state the following lemma, which will be useful in the error analysis. The proof of the following lemma follows from \cite[Lemma 6]{Hild_15_signorini}.


\begin{lemma}\label{lemma_req:para}
	Let $\sigma_e$ be the diameter of the two dimensional trace element $K_{\sigma} \cap\Gamma_C$, and $|S_{C}|>0$ and $|S_{NC}|>0$. Then the following estimations hold for $\mu_{n}$ and $\nabla \bar{q}$:
	\begin{align}
		\norm{\mu_n}_{0,K_{\sigma} \cap\Gamma_C}&\leq\frac{1}{|S_{NC}|^{1/2}}\;\sigma_e^{\tau-\frac{1}{2}}\; |\mu_n|_{\tau-\frac{3}{2},K_{\sigma} \cap\Gamma_C},\label{lemma_result_1:para}\\
		\norm{\mu_n}_{L^1(K_{\sigma} \cap\Gamma_C)}&\leq\frac{|S_{C}|^{1/2}}{|S_{NC}|^{1/2}}\;\sigma_e^{\tau-\frac{1}{2}}\; |\mu_n|_{\tau-\frac{3}{2},K_{\sigma} \cap\Gamma_C},\label{lemma_result_2:para}\\
		\norm{\nabla \bar{q}}_{0,K_{\sigma} \cap\Gamma_C}&\leq\frac{1}{|S_{C}|^{1/2}}\;\sigma_e^{\tau-\frac{1}{2}}\; |\nabla \bar{q}|_{\tau-\frac{3}{2},K_{\sigma} \cap\Gamma_C},\label{lemma_result_3:para}\\
		\norm{\nabla \bar{q}}_{L^1(K_{\sigma} \cap\Gamma_C)}&\leq\frac{|S_{NC}|^{1/2}}{|S_{C}|^{1/2}}\;\sigma_e^{\tau-\frac{1}{2}}\; |\nabla \bar{q}|_{\tau-\frac{3}{2},K_{\sigma} \cap\Gamma_C},\label{lemma_result_4:para}
	\end{align}
	where $\mu_n:=\lambda\frac{\partial \bar{q}}{\partial n}-\frac{\partial \bar{\phi}(\bar{q})}{\partial n}$ and $3/2<\tau\leq 2$. 
\end{lemma}
\begin{theorem}\label{Hilld_trick}
	For $\bar{q}\in H^{\tau}(\Omega \times I)$ with $3/2<\tau\leq 2$, it holds
	\begin{align*}
		|\lambda(\partial_t \bar{q},\partial_t (p_{\sigma}&-\bar{q}))_{I}+\lambda (\nabla \bar{q},\nabla (p_{\sigma}-\bar{q}))_I 
		-(\partial_t(p_{\sigma}-\bar{q}),\bar{\phi}(\bar{q}))_I\nonumber\\&-(\nabla(p_{\sigma}-\bar{q}),\nabla\bar{\phi}(\bar{q}))_I+(\bar{u}(\bar{q})-u_d,p_{\sigma}-\bar{q})_I|\leq C \sigma_e^{2(\tau-1)}.
	\end{align*}
	
\end{theorem}
\begin{proof}
	Integration by parts yields
	\begin{align}\label{Equn:9}
		\lambda(\partial_t \bar{q},\partial_t (p_{\sigma}-&\bar{q}))_{I}+\lambda (\nabla \bar{q},\nabla (p_{\sigma}-\bar{q}))_I 
		-(\partial_t(p_{\sigma}-\bar{q}),\bar{\phi}(\bar{q}))_I\nonumber\\&-(\nabla(p_{\sigma}-\bar{q}),\nabla\bar{\phi}(\bar{q}))_I+(\bar{u}(\bar{q})-u_d,p_{\sigma}-\bar{q})_I=\int_{\Gamma_C}\mu_n (p_{\sigma}-\bar{q}),
	\end{align} 
	where $\mu_n:=\rho\frac{\partial q}{\partial n}-\frac{\partial \phi}{\partial n}$.
	Now putting $p_{\sigma}=\mathcal{I}_{\sigma}\bar{q}+\mathcal{R}_{\sigma}\big(\mathcal{J}_{\sigma}(\bar{q}|_{\Gamma_C})-\mathcal{I}_{\sigma}(\bar{q}|_{\Gamma_C})\big)$ in \eqref{Equn:9} and using the property of $\mathcal{R}_{\sigma}$, we obtain
	\begin{align}\label{eq:contact_int:para}
		\int_{\Gamma_C}\mu_n (p_{\sigma}-\bar{q})=\int_{\Gamma_C}\mu_n (\mathcal{J}_{\sigma}\bar{q}-\bar{q})=\sum_{K_{\sigma}\in \mathcal{T}_{\sigma}} \int_{K_{\sigma} \cap\Gamma_C} \mu_n (\mathcal{J}_{\sigma}\bar{q}-\bar{q})ds.
	\end{align}	
	\noindent	
	Therefore it remains to estimate the following:
	\begin{equation}\label{eqn:estimate:para}
		\int_{K_{\sigma} \cap\Gamma_C} \mu_n (\mathcal{J}_{\sigma}\bar{q}-\bar{q})ds \quad \forall K_{\sigma}\in \mathcal{T}_{\sigma}.
	\end{equation}
	\noindent
	Let $K_{\sigma}$ be a fixed prism sharing a face with the boundary $\Gamma_C$ and  $\sigma_e$ be the diameter of the face $K_{\sigma} \cap\Gamma_C$ and obviously $|S_{C}|+|S_{NC}|=\textit{m}\sigma_e^2,$ where $\textit{m}$ is a fixed positive constant. Then, two cases can arise:
	\begin{itemize}
		\item[(a)] either $|S_{C}|$ or $|S_{NC}|$ equals zero,
		\item[(b)] both $|S_{C}|$ and $|S_{NC}|$  are positive.
	\end{itemize}
	
	It can be observed that the integral term in \eqref{eqn:estimate:para} vanishes for the first case (a). For the second case (b), we derive two estimations for the same error term \eqref{eqn:estimate:para}.
	
	\textit{The estimation of \eqref{eqn:estimate:para} related to $S_{NC}$:} A use of Cauchy--Schwarz inequality, estimation for \eqref{lemma_result_1:para} in Lemma \ref{lemma_req:para}, and standard estimation for the interpolation $\mathcal{J}_{\sigma}$ lead to
	\begin{align}
		\int_{K_{\sigma} \cap\Gamma_C} \mu_n (\mathcal{J}_{\sigma}\bar{q}-\bar{q})ds &\leq \norm{\mu_n}_{0,K_{\sigma} \cap\Gamma_C} \norm{\mathcal{J}_{\sigma}\bar{q}-\bar{q}}_{0,K_{\sigma} \cap\Gamma_C}\nonumber\\&\leq C \frac{1}{|S_{NC}|^\frac{1}{2}} \sigma_e^{\tau-\frac{1}{2}} |\mu_n|_{\tau-\frac{3}{2},K_{\sigma} \cap\Gamma_C} \sigma_e^{\tau-\frac{1}{2}} |\nabla \bar{q}|_{\tau-\frac{3}{2},K_{\sigma} \cap\Gamma_C}\nonumber\\& \leq  C \frac{1}{|S_{NC}|^\frac{1}{2}} \sigma_e^{2(\tau-\frac{1}{2})} \big(|\mu_n|^2_{\tau-\frac{3}{2},K_{\sigma} \cap\Gamma_C} + |\nabla \bar{q}|^2_{\tau-\frac{3}{2},K_{\sigma} \cap\Gamma_C}\big).\label{para:1}
	\end{align}
	\textit{Estimation for \eqref{eqn:estimate:para} related to $S_{C}$:} Using the  estimation for $\mathcal{J}_{\sigma}$ in \eqref{Estm:12} and estimations \eqref{lemma_result_1:para} and \eqref{lemma_result_4:para}, we obtain
	\begin{align}
		\int_{K_{\sigma} \cap\Gamma_C} \mu_n (\mathcal{J}_{\sigma}\bar{q}&-\bar{q})ds \leq \norm{\mu_n}_{0,K_{\sigma} \cap\Gamma_C} \norm{\mathcal{J}_{\sigma}\bar{q}-\bar{q}}_{0,K_{\sigma} \cap\Gamma_C}\nonumber\\&\leq C \frac{1}{|S_{NC}|^\frac{1}{2}} \sigma_e^{\tau-\frac{1}{2}} |\mu_n|_{\tau-\frac{3}{2},K_{\sigma} \cap\Gamma_C} \norm{\nabla \bar{q}}_{L^1(K_{\sigma} \cap\Gamma_C)}\nonumber\\&
		\leq C \frac{1}{|S_{NC}|^\frac{1}{2}} \sigma_e^{\tau-\frac{1}{2}} |\mu_n|_{\tau-\frac{3}{2},K_{\sigma} \cap\Gamma_C}\frac{|S_{NC}|^{1/2}}{|S_{C}|^{1/2}}\;\sigma_e^{\tau-\frac{1}{2}}\; |\nabla \bar{q}|_{\tau-\frac{3}{2},K_{\sigma} \cap\Gamma_C}\nonumber\\&
		\leq  C \frac{1}{|S_{C}|^\frac{1}{2}} \sigma_e^{2(\tau-\frac{1}{2})} \big(|\mu_n|^2_{\tau-\frac{3}{2},K_{\sigma} \cap\Gamma_C} + |\nabla \bar{q}|^2_{\tau-\frac{3}{2},K_{\sigma} \cap\Gamma_C}\big).\label{para:2}
	\end{align}
	It is easy to observe that either $|S_{NC}|$ or $|S_{C}|$ is greater than or equal to $\textit{m}\sigma_e^2/2$. Then, choosing
	the appropriate estimation \eqref{para:1} or \eqref{para:2}, we obtain
	\begin{equation*}
		\int_{K_{\sigma} \cap\Gamma_C} \mu_n (\mathcal{J}_{\sigma}\bar{q}-\bar{q})ds\leq  C  \sigma_e^{2(\tau-1)} \big(|\mu_n|^2_{\tau-\frac{3}{2},K_{\sigma} \cap\Gamma_C} + |\nabla \bar{q}|^2_{\tau-\frac{3}{2},K_{\sigma} \cap\Gamma_C}\big).
	\end{equation*}
	Summing over all $K_{\sigma}$ sharing a face with $\Gamma_{C}$ and applying the trace theorem, we get
	\begin{equation*}
		\int_{\Gamma_C} \mu_n (\mathcal{J}_{\sigma}\bar{q}-\bar{q})ds\leq  C \sigma_e^{2(\tau-1)} \big(|\mu_n|^2_{\tau-\frac{3}{2},\Gamma_C} + |\nabla \bar{q}|^2_{\tau-\frac{3}{2},\Gamma_C}\big)\leq
		\sigma_e^{2(\tau-1)} \norm{\bar{q}}^2_{\tau,\Omega\times I}.
	\end{equation*}
	This completes the proof.
\end{proof}

In the following theorem, we derive the energy error estimate for the control and $L^2$-error estimate of the state variable.

\begin{theorem} [Error estimate of control variable]\label{control_estimate} There holds
	\begin{align*}
		\lambda |\bar{q}-\bar{q}_{\sigma}|_{1,\Omega\times I}+&\norm{\bar{u}(\bar{q})-\bar{u}_{kh}(\bar{q}_{\sigma})}_I\leq C (\sigma_e^{\tau-1}+h^{\tau-1}+k).
	\end{align*}
	
\end{theorem}
\begin{proof} Recall the result of Lemma \ref{Estm:Control}:
	\begin{align}
		\lambda |\bar{q}-\bar{q}_{\sigma}|^2_{1,\Omega\times I}+&\norm{\bar{u}(\bar{q})-\bar{u}_{kh}(\bar{q}_{\sigma})}^2_I
		\leq 
		\big[\lambda(\partial_t \bar{q},\partial_t (p_{\sigma}-\bar{q}))_{I}+\lambda (\nabla \bar{q},\nabla (p_{\sigma}-\bar{q}))_I\nonumber\\& 
		-(\partial_t(p_{\sigma}-\bar{q}),\bar{\phi}(\bar{q}))_I-(\nabla(p_{\sigma}-\bar{q}),\nabla\bar{\phi}(\bar{q}))_I+(\bar{u}(\bar{q})-u_d,p_{\sigma}-\bar{q})_I\big]\nonumber\\&+\norm{\bar{q}-\bar{q}_{\sigma}}^2_{1,\Omega\times I}
		+\norm{\nabla(\phi_{kh}(\bar{q})-\bar{\phi}(\bar{q}))}^2_I+\norm{w_{kh}(\bar{q})-\bar{w}(\bar{q})}^2_I\label{L2:estimate}
	\end{align}
	for all $p_{\sigma}\in Q^{\sigma}_{ad}$. From Theorem \ref{Hilld_trick}, we obtain an estimation for the first term of the above equation \eqref{L2:estimate} as
	\begin{align}\label{Hilld1_trem}
		|\lambda(\partial_t \bar{q},\partial_t (p_{\sigma}-\bar{q}))_{I}&+\lambda (\nabla \bar{q},\nabla (p_{\sigma}-\bar{q}))_I 
		-(\partial_t(p_{\sigma}-\bar{q}),\bar{\phi}(\bar{q}))_I\nonumber\\&-(\nabla(p_{\sigma}-\bar{q}),\nabla\bar{\phi}(\bar{q}))_I+(\bar{u}(\bar{q})-u_d,p_{\sigma}-\bar{q})_I|\leq C\sigma_e^{2(\tau-1)}.
	\end{align}
	For the second term in the right hand side of \eqref{L2:estimate},
	we take $p_{\sigma}=\mathcal{I}_{\sigma}\bar{q}+\mathcal{R}_{\sigma}\big(\mathcal{J}_{\sigma}(\bar{q}|_{\Gamma_C})-\mathcal{I}_{\sigma}(\bar{q}|_{\Gamma_C})\big)$.
	The continuity of the extension operator $\mathcal{R}_{\sigma}$ and an inverse inequality yield
	\begin{align}\label{Best_approx}
		\norm{\bar{q}-p_{\sigma}}_{1,\Omega\times I}&\leq \norm{\bar{q}-\mathcal{I}_{\sigma}\bar{q}}_{1,\Omega\times I}+\norm{\mathcal{R}_{\sigma}\big(\mathcal{J}_{\sigma}(\bar{q}|_{\Gamma_C})-\mathcal{I}_{\sigma}(\bar{q}|_{\Gamma_C})\big)}_{1,\Omega\times I}\nonumber\\
		&\leq C \sigma_e^{\tau-1}\norm{\bar{q}}_{\tau,\Omega\times I}+C \norm{\mathcal{J}_{\sigma}(\bar{q}|_{\Gamma_C})-\mathcal{I}_{\sigma}(\bar{q}|_{\Gamma_C})}_{\frac{1}{2},\Gamma_C}\nonumber\\
		&\leq C \sigma_e^{\tau-1}\norm{\bar{q}}_{\tau,\Omega\times I}+C\sigma_e^{-\frac{1}{2}} \norm{\mathcal{J}_{\sigma}(\bar{q}|_{\Gamma_C})-\mathcal{I}_{\sigma}(\bar{q}|_{\Gamma_C})}_{0,\Gamma_C}\nonumber\\
		&\leq C \sigma_e^{\tau-1}\norm{\bar{q}}_{\tau,\Omega\times I}+C\sigma_e^{-\frac{1}{2}} \norm{\bar{q}|_{\Gamma_C}-\mathcal{J}_{\sigma}(\bar{q}|_{\Gamma_C})}_{0,\Gamma_C}\nonumber\\&\;\;+C\sigma_e^{-\frac{1}{2}}\norm{\bar{q}|_{\Gamma_C}-\mathcal{I}_{\sigma}(\bar{q}|_{\Gamma_C})}_{0,\Gamma_C}\nonumber\\
		&\leq C\sigma_e^{\tau-1}\norm{\bar{q}}_{\tau,\Omega\times I}.
	\end{align}
	The above estimations \eqref{Hilld1_trem}, \eqref{Best_approx} and Theorems \ref{estm_adj_fixed_q} \& \ref{Error_estm_state_fixed_q} lead to the required result.
\end{proof}

\begin{theorem}[Error estimate of state variable]\label{state_estimate:para} There holds,
	\begin{align*}
		\norm{\nabla(\bar{u}(\bar{q})-\bar{u}_{kh}(\bar{q}_\sigma))}_I&\leq C (\sigma_e^{\tau-1}+h^{\tau-1}+k).
	\end{align*}
\end{theorem}
\begin{proof}
	The triangle inequality gives
	\begin{align}\label{Estm:19}
		\norm{\nabla(\bar{u}(\bar{q})-\bar{u}_{kh}(\bar{q}_\sigma))}_I
		&\leq \norm{\nabla(\bar{u}(\bar{q})-u_{kh}(\bar{q}))}_I+\norm{\nabla(u_{kh}(\bar{q})-\bar{u}_{kh}(\bar{q}_\sigma))}_I.
	\end{align}
	For the first term of \eqref{Estm:19}, we use the splitting $\bar{u}(\bar{q})=\bar{w}(\bar{q})+\bar{q}$ from the equation \eqref{eq:state_0a:para} and $u_{kh}(\bar{q})=w_{kh}(\bar{q})+\bar{q}$ from the equation \eqref{Eq:fully_disc_state} to obtain
	\begin{align}\label{Estm:21}
		\norm{\nabla(\bar{u}(\bar{q})-u_{kh}(\bar{q}))}_I=\norm{\nabla(\bar{w}(\bar{q})-w_{kh}(\bar{q}))}_I.	
	\end{align}
	For the second term of \eqref{Estm:19}, we use the splitting $u_{kh}(\bar{q})=w_{kh}(\bar{q})+\bar{q}$ from the equation \eqref{Eq:fully_disc_state} and $\bar{u}_{kh}(\bar{q}_\sigma)=\bar{w}_{kh}(\bar{q}_\sigma)+\bar{q}_\sigma$ from the equation \eqref{full_disc_state}. Hence, we have
	\begin{align*}
		u_{kh}(\bar{q})-\bar{u}_{kh}(\bar{q}_\sigma)=w_{kh}(\bar{q})-\bar{w}_{kh}(\bar{q}_\sigma)+\bar{q}-\bar{q}_\sigma.	
	\end{align*}
	Using the stability estimate \eqref{stability_estimate_state} of the fully discrete state equation, we obtain 
	\begin{align}\label{Estm:20}
		\norm{\nabla(u_{kh}(\bar{q})-\bar{u}_{kh}(\bar{q}_\sigma))}_I&\leq \norm{\nabla(w_{kh}(\bar{q})-\bar{w}_{kh}(\bar{q}_\sigma))}_I+\norm{\nabla(\bar{q}-\bar{q}_\sigma)}_I\nonumber\\
		&\leq C |\bar{q}-\bar{q}_\sigma|_{1,\Omega\times I}.	
	\end{align}
	Putting \eqref{Estm:21} and \eqref{Estm:20} in \eqref{Estm:19}, we get
	\begin{align}\label{Estm:22}
		\norm{\nabla(\bar{u}(\bar{q})-\bar{u}_{kh}(\bar{q}_\sigma))}_I
		&\leq \norm{\nabla(\bar{w}(\bar{q})-w_{kh}(\bar{q}))}_I+|\bar{q}-\bar{q}_\sigma|_{1,\Omega\times I}.
	\end{align}
	The estimations for Theorem \ref{Error_estm_state_fixed_q} and Theorem \ref{control_estimate} lead to the required result.
\end{proof}

\begin{theorem} [Error estimate of adjoint state]\label{adj_state_estm:para} Let $\bar{\phi}(\bar{q})$ be the solution of \eqref{eq:adjstate_a:para} and $\bar{\phi}_{kh}(\bar{q}_\sigma)$ be the solution of \eqref{full_disc_adj_state}. Then there holds,
	\begin{align*}
		\norm{\nabla(\bar{\phi}(\bar{q})-\bar{\phi}_{kh}(\bar{q}_\sigma))}_I\leq C ( k+ h^{\tau-1}+\sigma_e^{\tau-1}).
	\end{align*}
\end{theorem}
\begin{proof}
	Introducing the auxiliary solution $\phi_{kh}(\bar{q})$ satisfying \eqref{Aux_adj_state}, we obtain
	\begin{align}\label{Estm:28}
		\norm{\nabla(\bar{\phi}(\bar{q})-\bar{\phi}_{kh}(\bar{q}_\sigma))}_I
		&\leq \norm{\nabla(\bar{\phi}(\bar{q}) -\phi_{kh}(\bar{q}))}_I+\norm{\nabla(\phi_{kh}(\bar{q})-\bar{\phi}_{kh}(\bar{q}_\sigma))}_I.
	\end{align}
	For the second term of \eqref{Estm:28}, we use the stability estimate of the discrete adjoint solution \eqref{stability_estimate_costate} to obtain
	\begin{align}\label{Estm:29}
		\norm{\nabla(\phi_{kh}(\bar{q})-\bar{\phi}_{kh}(\bar{q}_\sigma))}_I&\leq C \norm{\bar{u}(\bar{q})-\bar{u}_{kh}(\bar{q}_\sigma)}_I.	
	\end{align}
	The estimations for Theorem \ref{estm_adj_fixed_q} and Theorem \ref{control_estimate} lead to the required result.
\end{proof}
\begin{remark}\label{estm_less_reg:para}
	Note that the optimal control $\bar{q}$ satisfies a simplified Signorini problem. The regularity of the solution of Signorini problem gets impaired due to many reasons, for example regularity of the data, the mixed boundary conditions
	(e.g., Neumann-Dirichlet transitions), the corners in polygonal domains and the Signorini condition which generates singularities at contact-noncontact transition points which we have discussed in the Remark \ref{Rem_reg_signorini:para}. So, there is a possibility that the solution could be less regular i.e., $\bar{q}\in H^{\tau}(\Omega \times I),$ where $ 1<\tau\leq3/2$. Then all the above \textit{a priori} estimates hold true except the Theorem \ref{Hilld_trick}. It is clear that if the solutions have the above regularity then \eqref{Equn:9} is not true because the right hand side of \eqref{Equn:9} does not make sense. So, to estimate the term 
	\begin{align}
		\lambda(\partial_t \bar{q},\partial_t (p_{\sigma}-\bar{q}))_{I}+&\lambda (\nabla \bar{q},\nabla (p_{\sigma}-\bar{q}))_I 
		-(\partial_t(p_{\sigma}-\bar{q}),\bar{\phi}(\bar{q}))_I\nonumber\\&-(\nabla(p_{\sigma}-\bar{q}),\nabla\bar{\phi}(\bar{q}))_I+(\bar{u}(\bar{q})-u_d,p_{\sigma}-\bar{q})_I,
	\end{align} 
	we use the following idea:
	\begin{align}\label{Estm:30_1}
		\lambda(\partial_t \bar{q},&\partial_t (p_{\sigma}-\bar{q}))_{I}+\lambda (\nabla \bar{q},\nabla (p_{\sigma}-\bar{q}))_I 
		-(\partial_t(p_{\sigma}-\bar{q}),\bar{\phi}(\bar{q}))_I\nonumber\\&-(\nabla(p_{\sigma}-\bar{q}),\nabla\bar{\phi}(\bar{q}))_I+(\bar{u}(\bar{q})-u_d,p_{\sigma}-\bar{q})_I=\langle\mu_n, (p_{\sigma}-\bar{q})\rangle_{\epsilon,\Gamma_C}\nonumber\\& \leq \norm{\mu_n}_{H^{\epsilon}(\Gamma_C)'} \norm{p_{\sigma}-\bar{q}}_{\epsilon,\Gamma_C},
	\end{align}
	where $\epsilon=3/2-\tau$ and $H^{\epsilon}(\Gamma_C)'$ denotes the dual of $H^{\epsilon}(\Gamma_C)$ (see \cite{Belhachmi:2003}).	
	Choosing $p_{\sigma}=\mathcal{I}_{\sigma}\bar{q},$ we have 
	$\norm{\bar{q}-\mathcal{I}_{\sigma}\bar{q}}_{\epsilon,\Gamma_C}\leq C \sigma^{2\tau-2} \norm{\bar{q}}_{\tau,\Omega}$. Using the trace estimate (discussed in Section \ref{sec:ModelProblem:para}), we have $\norm{\mu_n}_{H^{\epsilon}(\Gamma_C)'}\leq C \norm{\bar{q}}_{\tau,\Omega}$. Putting all these estimates in \eqref{Estm:30_1}, we have
	\begin{align}
		|\lambda(\partial_t \bar{q},&\partial_t (p_{\sigma}-\bar{q}))_{I}+\lambda (\nabla \bar{q},\nabla (p_{\sigma}-\bar{q}))_I 
		-(\partial_t(p_{\sigma}-\bar{q}),\bar{\phi}(\bar{q}))_I\nonumber\\&-(\nabla(p_{\sigma}-\bar{q}),\nabla\bar{\phi}(\bar{q}))_I+(\bar{u}(\bar{q})-u_d,p_{\sigma}-\bar{q})_I| \leq C \sigma^{2\tau-2} \norm{\bar{q}}^2_{\tau,\Omega}.
	\end{align}
	Thus, we have an optimal order (up to the regularity) of convergence of the term \eqref{Equn:9}. Hence, all the error estimations (control, state and adjoint state) show the optimal order of convergence (up to the regularity of the solutions). 	  
\end{remark}
So, it is clear from the Remark \ref{estm_less_reg:para} that our error analysis also works for the solutions with low regularity.

\section{Numerical Experiments}\label{sec5:para}
In this section, we validate the a priori error estimates for the
error in state, adjoint state and control variables numerically.  We use primal-dual active set strategy (see \cite{trolzstch:2005:Book}) in combination with conjugate gradient method (see, \cite{Meid_Vexler_08_II,Meid_Vexler_08_I}) to solve the optimal control problem. For the computations we construct a model problem with known solutions. In order to accomplish this, we consider the following cost functional $\tilde{J}$ defined by 
$$\tilde{J}(u,q):= \frac{1}{2}\norm{u-u_d}_I^2 + \frac{\lambda}{2}|q-q_d|^2_{1,\Omega\times I}, \quad
w\in Q,\;p\in Q_{ad},$$

\noindent
for some given function $q_d$. Then the minimization problem reads: Find $(u,q)\in Q \times Q_{ad}$ such that
$$ \tilde{J}(\bar{u}(\bar{q}),\bar{q})=\min_{(u,q)\in Q\times Q_{ad}}  \tilde{J}(u,q)$$
subject to the condition that $(u,q)\in (X+Q) \times Q_{ad}$  satisfies the state equation \eqref{Eq:state_weak_form}.  Then the discrete optimality system finds
$(\bar{u}_{kh}(\bar{q}_{\sigma}),\bar{\phi}_{kh}(\bar{q}_{\sigma}),\bar{q}_{\sigma})\in (X^{0,1}_{k,h}+Q_{\sigma})\times X_{k,h}^{0,1} \times Q^{\sigma}_{ad}$ such that

\begin{subequations}
	\begin{align}
		\bar{u}_{kh}(\bar{q}_{\sigma})=&\bar{w}_{kh}(\bar{q}_{\sigma})+\bar{q}_{\sigma}\quad \bar{w}_{kh}(\bar{q}_{\sigma}) \in X_{k,h}^{0,1}\\
		B(\bar{w}_{kh}(\bar{q}_{\sigma}),v_{kh})&=(f,v_{kh})_I+(u_0,v^{+}_{kh,0})-B(\bar{q}_{\sigma},v_{kh})\quad \forall v_{kh} \in X^{0,1}_{k,h}\\
		B(v_{kh},\bar{\phi}_{kh}(\bar{q}_{\sigma}),)&=\big(\bar{u}_{kh}(\bar{q}_{\sigma})-u_d,v_{kh}\big)_I\quad \forall v_{kh} \in X^{0,1}_{k,h}\\
		\lambda(\partial_t \bar{q}_{\sigma},\partial_t (p_{\sigma}-\bar{q}_{\sigma}))_{I}&+\lambda (\nabla \bar{q}_{\sigma},\nabla(p_{\sigma}-\bar{q}_{\sigma}))_I\geq  (\partial_t(p_{\sigma}-\bar{q}_{\sigma}), \bar{\phi}_{kh}(\bar{q}_{\sigma}))_I\nonumber\\&+(\nabla \bar{\phi}_{kh}(\bar{q}_{\sigma}),\nabla(p_{\sigma}-\bar{q}_{\sigma}))_I-(\bar{u}_{kh}(\bar{q}_{\sigma})-u_d,p_{\sigma}-\bar{q}_{\sigma})_I\nonumber\\&+\lambda(\partial_t q_d,\partial_t (p_{\sigma}-\bar{q}_{\sigma}))_{I}+\lambda (\nabla q_d , \nabla (p_{\sigma}-\bar{q}_{\sigma}))_I,
	\end{align}	
\end{subequations}
for all $p_{\sigma} \in Q^{\sigma}_{ad}$

\begin{example}\label{Ex.1.1} 
	Let the computational domain be $\Omega \times I:=(0,1)^2\times(0,1)$, $\Gamma_C:=\gamma_C \times (0,1)$, and $\Gamma_D:=\partial(\Omega \times I) \setminus \Gamma_C$ where $\gamma_C:= (0,1)\times \{0\}$. We choose the exact solutions as follows:
	\begin{align*}
		u(x,y)&=x\exp{(y)}\;(1-x)(1-y) t (1-t),\\
		\phi(x,y)&=(x^2-x^3)(y^2-y^3)t(1-t),\\
		q(x,y)&=x\exp{(y)}\;(1-x)(1-y) t (1-t),
	\end{align*}
	and set the data as
	\begin{align*}
		f&=\partial _t u-\Delta u,\\
		u_d&=u+\partial _t \phi+\Delta\phi,\\
		q_d&=q,\\
		\lambda&=10^{-3}, q_a=0, q_b=0.8.
	\end{align*}
\end{example}
In this numerical experiments, we consider a sequence of uniformly refined meshes. The spatial domain $\Omega$ is subdivided by regular triangular elements and the time interval is partitioned by equally spaced time steps. To discretize the state and adjoint state we use piecewise linear and continuous finite elements for  spatial discretization and piecewise constant elements for temporal discretization. For the discretization of control we use linear prismatic Lagrange finite elements. We compute the errors in state, adjoint state, and control on the above mentioned uniformly refined meshes.
The empirical convergence rate is defined by
\begin{equation*}
	\texttt{rate}(\ell):= \log(e_{\ell}/e_{\ell-1})/\log(\mu_\ell/\mu_{\ell-1}), \quad \text{for} \; \ell=1,2,3,...
\end{equation*} 
where $e_{\ell}$ and $\mu_\ell$ denote respectively the error and the discretization parameter at $\ell$-th level. Let $N$ denote the number of sub-intervals for the time interval $\bar{I}$. In Table \ref{table1.1_a}, we have shown the rate of convergence of state and adjoint state in the energy norm with respect to the space parameter $h$. Table \ref{table1.1_b} shows the rate of convergence of state and adjoint state in the $L^2$-norm with respect to the time parameter $k$. In Table \ref{table1.1_c}, we have shown rate of convergence of the control variable in the energy norm with respect to the control discretization parameter $\sigma:=\sqrt{h^2+k^2}$.

\begin{table}[h!]
	\begin{center}
		\footnotesize
		\vspace{0.2cm}
		\caption{Errors and  rates of convergence of state and adjoint state w.r.t. space parameter $h$ for Example \ref{Ex.1.1}.}
		\label{table1.1_a}
		\begin{tabular}{|c|c|c|c|c|c|}
			\hline
			$N$&	$h$ &$\norm{\nabla(\bar{u}(\bar{q})-\bar{u}_{kh}(\bar{q}_{\sigma}))}_I$&\texttt{rate} & $\norm{\nabla(\bar{\phi}(\bar{q})-\bar{\phi}_{kh}(\bar{q}_{\sigma}))}_I$ & \texttt{rate} \\
			\hline                     
			4&0.2500  & 0.02610199  &  --------  & 0.00690363  & ------\\
			6&0.1250  & 0.01401513 &   0.8971    & 0.00341977  & 1.0134\\
			12&0.0625 & 0.00707057 &    0.9870   & 0.00165730  & 1.0450\\
			23&0.0312 & 0.00357310 &   0.9846    & 0.00081186   & 1.0295\\
			46&0.0156 & 0.00178706 &  0.9995     & 0.00040030   & 1.0201\\
			\hline
		\end{tabular}
	\end{center}
\end{table}

\begin{table}[h!]
	\begin{center}
		\footnotesize
		\vspace{0.2cm}
		\caption{Errors and  rates of convergence of state and adjoint state w.r.t. time parameter $k$ for Example \ref{Ex.1.1}.}
		\label{table1.1_b}
		\begin{tabular}{|c|c|c|c|c|c|}
			\hline
			$N$&	$h$ &$\norm{\nabla(\bar{u}(\bar{q})-\bar{u}_{kh}(\bar{q}_{\sigma}))}_I$& \texttt{rate}   &$\norm{\nabla(\bar{\phi}(\bar{q})-\bar{\phi}_{kh}(\bar{q}_{\sigma}))}_I$  & \texttt{rate} \\
			
			\hline

			4&0.2500  & 0.02610199  & -------     &0.00690363 & ------\\
			6&0.1250  & 0.01401513 &   1.5337     &0.00341977   & 1.7325 \\
			12&0.0625 & 0.00707057 &   0.9870     &0.00165730  & 1.0450 \\
			23&0.0312 & 0.00357310 &   0.9870     &0.00081186    & 1.0968\\
			46&0.0156 & 0.00178706 &   0.9995     &0.00040030  & 1.0968\\
			
			\hline
		\end{tabular}
	\end{center}
\end{table}

\begin{table}[h!]
	\begin{center}
		\footnotesize
		\vspace{0.2cm}
		\caption{Errors and  rates of convergence of control variable for Example \ref{Ex.1.1}.}
		\label{table1.1_c}
		\begin{tabular}{|c|c|c|c|}
			\hline
			$N$&	$\sigma$ &$|\bar{q}-\bar{q}_{\sigma}|_{1,\Omega\times I}$& \texttt{rate} \\
			
			\hline
			
			4& 0.3535  &0.09476646  & ------ \\
			6& 0.2083  &0.05263751  & 1.1117\\
			12&0.1041  &0.02631136 & 1.0004\\
			23& 0.0535 &0.01342729  & 1.0004\\
			46& 0.0267 &0.00671436  &0.9998\\
			
			\hline
		\end{tabular}
	\end{center}
\end{table}

\section*{Conclusions}\label{sec6:para}
We address the energy approach to solve the Dirichlet boundary control problem governed by the linear parabolic equation. Since we have chosen the control from a closed convex subset of  $H^1(\Omega\times (0,T))$, the optimal control satisfies a simplified Signorini problem in three dimensional domain $\Omega\times (0,T)$. For the discretization, we use conforming prismatic Lagrange finite elements for the control. We derive the optimal order of convergence for the error in control, state, and adjoint state. Our numerical experiments confirm the theoretical results.


\bibliographystyle{amsplain}
\bibliography{references}
\end{document}